\documentclass[11pt]{article}

\usepackage{tikz-cd}
\usepackage[a4paper,margin=1in]{geometry}
\usepackage{lmodern}
\usepackage{microtype}
\usepackage{amsmath,amssymb,amsthm,mathtools,stmaryrd,mathrsfs}
\usepackage[backend=bibtex,style=alphabetic,maxnames=10,minnames=10,maxbibnames=99,minalphanames=4,maxalphanames=10]{biblatex}
\usepackage[colorlinks=true,linkcolor=blue!55!black,citecolor=blue!55!black,urlcolor=blue!55!black]{hyperref}
\usepackage[nameinlink]{cleveref}

\newtheorem{theorem}{Theorem}[section]
\newtheorem{proposition}[theorem]{Proposition}
\newtheorem{lemma}[theorem]{Lemma}
\newtheorem{corollary}[theorem]{Corollary}
\newtheorem*{mainresult}{Theorem}

\theoremstyle{definition}
\newtheorem{definition}[theorem]{Definition}
\newtheorem{hypothesis}[theorem]{Hypothesis}

\newtheorem{remark}[theorem]{Remark}
\newtheorem{example}[theorem]{Example}

\newcommand{\remarkendmark}{\ensuremath{\diamond}}
\let\proofendmark\qedsymbol
\AddToHook{env/remark/begin}{%
  \renewcommand{\qedsymbol}{\remarkendmark}\pushQED{\qed}}
\AddToHook{env/remark/end}{\popQED}
\AddToHook{env/proof/begin}{\let\qedsymbol\proofendmark}

\newcommand{\F}{\mathbb F}
\newcommand{\cO}{\mathcal O}

\newcommand{\mD}{\mathfrak m_D}

\newcommand{\gr}{\operatorname{gr}}
\newcommand{\Fil}{\operatorname{Fil}}
\newcommand{\rank}{\operatorname{rank}}
\newcommand{\im}{\operatorname{Im}}
\newcommand{\Hom}{\operatorname{Hom}}
\newcommand{\GL}{\operatorname{GL}}
\newcommand{\Nrd}{\operatorname{Nrd}}

\newcommand{\JL}{\operatorname{JL}}
\newcommand{\Tau}{\mathcal T}
\newcommand{\CC}{\operatorname{CC}}
\newcommand{\Supp}{\operatorname{Supp}}

\newcommand{\Min}{\operatorname{Min}}
\newcommand{\length}{\operatorname{length}}

\title{On the Irreducibility of Locally Analytic Vectors in Certain Banach Representations of \(D^\times\)}
\author{Yu Xiao}
\date{\today}

\begin{document}
\maketitle

\begin{abstract}
      We investigate the structure of locally analytic representations of $p$-adic groups, 
      and establish a sufficient condition for the irreducibility of locally analytic vectors in certain Banach representations of \(D^\times\), where $D$ is the quaternion division algebra over $\mathbb{Q}_p$.
      As an application, we prove that \(\JL(\rho)^{\mathrm{an}}/\JL(\rho)^{1\text{-}\mathrm{alg}}\) is a topologically irreducible locally analytic representation of $D^\times$
      under certain genericity conditions on \(\rho\). Here \(\JL(\rho)\) is the \(\rho(-1)\)-multiplicity space in the image of \(\Pi(\rho)\) under Scholze's functor,
      and \(\Pi(\rho)\) is the representation associated to $\rho$ by the \(p\)-adic local Langlands correspondence.
\end{abstract}
\begingroup
\hypersetup{linkcolor=black}
\tableofcontents
\endgroup

\section*{Introduction}
\addcontentsline{toc}{section}{Introduction}

Let $p\geq 5$, and let $E/\mathbb Q_p$ be a sufficiently large finite
extension. We write $\cO_E$, $\varpi_E$, and $k$ for its ring of integers,
uniformizer, and residue field.  We assume $\mathbb F_{p^2}\subset k$.

Let $D$ be the nonsplit quaternion algebra over $\mathbb Q_p$, with
maximal order $\cO_D$ and a uniformizer $\varpi_D$ satisfying
$\varpi_D^2=p$.  Put
\[
 U_D^1=1+\varpi_D\mathcal{O}_D,
 \qquad Z_D^1=U_D^1\cap \mathbb Q_p^\times,
 \qquad
 H=[\mathbb F_{p^2}^{\times}]\subset \cO_D^\times,
\]
where the brackets denote the Teichm\"uller lift, so that $H$ maps
isomorphically onto $\cO_D^\times/U_D^1$.
All representations below have a fixed unitary central character.
We write $\zeta:Z(D^\times)\to\cO_E^\times$ for the character on the
module under consideration; passage to its contragredient replaces
$\zeta$ by $\zeta^{-1}$.

Let
\(
 \rho:\mathrm{Gal}_{\mathbb Q_p}\longrightarrow\GL_2(E)
\)
be a two-dimensional absolutely irreducible representation.  Fix a
$\mathrm{Gal}_{\mathbb Q_p}$-stable lattice $T_\rho\subset\rho$, and write
$\bar\rho=T_\rho/\varpi_ET_\rho$.  Let $\Pi(\rho)$ be the admissible
unitary Banach representation of $\GL_2(\mathbb Q_p)$ attached to $\rho$
by Colmez's $p$-adic local Langlands correspondence \cite{Colmez2010}.
We assume that the block $\mathscr B(\bar\rho)$ attached to $\bar\rho$ is
of type~(I) or type~(II) and is generic in the sense of Hu--Wang.  

Let $\pi(\bar\rho)$ be the representation attached to $\bar\rho$ by the
mod-$p$ local Langlands correspondence, and let $\mathcal S^i$ denote
Scholze's functors \cite{Scholze2018}. We use Hu--Wang's multiplicity
space $\JL(\bar\rho)$, characterized by
\[
 \mathcal S^1(\pi(\bar\rho))
 \cong\bar\rho(-1)\otimes_k\JL(\bar\rho)
\]
by \cite[Theorem~4.20]{HuWangPadic}. 
For the above generic blocks, \cite[Theorem~4.19]{HuWangPadic} shows that $\gr_{\mD}\bigl(\JL(\bar\rho)^\vee\bigr)$ is of the form
\[
\bigoplus_{\chi\in W_D(\bar{\rho})}
 \chi^\vee\otimes_k A_0/\mathfrak{a}(\chi),
 \qquad A_0=k[y,z]/(yz),
\]
as a module over the central reduction 
$\gr_{\mD}^*k\llbracket U_D^1/Z_D^1\rrbracket$, where the ideals
$\mathfrak a(\chi)$ are those of \cite[Definition~4.18]{HuWangPadic}.
Its $D$-multiplicity is 4.

Choose a $\GL_2(\mathbb Q_p)$-stable open bounded lattice
$\Theta_\rho\subset\Pi(\rho)$.  Using compatible integral models for
Scholze's degree-one functor, set \[\check{\mathcal{S}}^1(\Pi(\rho))\coloneqq\bigl(\varprojlim_n\mathcal{S}^1(\Theta_\rho/\varpi_E^n)\bigr)[1/\varpi_E],\quad\JL(\rho)=\mathrm{Hom}_{E[\mathrm{Gal}_{\mathbb{Q}_p}]}(\rho(-1), \check{\mathcal{S}}^1(\Pi(\rho))). \]

Write $\JL(\rho)^{1\text{-}\mathrm{alg}}$ for the vectors locally
algebraic for $D^{\times,1}:=\ker(\Nrd:D^\times\to\mathbb Q_p^\times)$.
By \cite[Theorem~4.23]{LiuWang2025},
\[
 \JL(\rho)^{\mathrm{lalg}}
 \subseteq\JL(\rho)^{1\text{-}\mathrm{alg}},
 \qquad
 \dim_E\JL(\rho)^{1\text{-}\mathrm{alg}}<\infty.
\]
If $\rho$ is Hodge--Tate with distinct integral weights $a>b$, the
infinitesimal-character calculation in
\cite[\S8.3 and Corollary~8.16]{DPS23} implies
\[
 \JL(\rho)^{\mathrm{lalg}}
 =\JL(\rho)^{1\text{-}\mathrm{alg}}.
\]

For arbitrary $\rho$, we set
\[
 \Tau=\Tau(\rho):=\JL(\rho)/\JL(\rho)^{1\text{-}\mathrm{alg}}.
\]
The denominator is finite-dimensional by \cite[Theorem~5.3]{LiuWang2025}, hence closed. 
By \cite[Corollary~4.20 and Theorem~6.4]{LiuWang2025}, the Banach
representations $\JL(\rho)$ and $\Tau$ are topologically of finite
length. In this article we focus on the locally analytic vectors of $\Tau$:

 \[\Tau^{\mathrm{an}}
 \cong
 \JL(\rho)^{\mathrm{an}}/
 \JL(\rho)^{1\text{-}\mathrm{alg}}.\]

If the difference of the
Hodge--Tate--Sen weights of $\rho$ is not a nonzero integer, then by \cite[Theorem~8.17]{DPS23},
\[
 \JL(\rho)^{1\text{-}\mathrm{alg}}
 =\JL(\rho)^{\mathrm{lalg}}=0
\] 
On the other hand, suppose $\rho$ is the local
representation in \cite[Theorem~3.20]{LiSuWu2026} and its associated
smooth representation $\pi_p$ of $\GL_2(\mathbb Q_p)$ is an irreducible
principal series. According to the classical-cohomology identification recalled
before \cite[Theorem~3.14]{LiSuWu2026} and the eigenspace
decomposition \cite[(56)]{LiSuWu2026}, together with classical
Jacquet--Langlands, we deduce that
\[
 \JL_{\mathrm{cl}}(\pi_p)=0
 \quad\Longrightarrow\quad
 \JL(\rho)^{\mathrm{lalg}}
 =\JL(\rho)^{1\text{-}\mathrm{alg}}=0.
\]
The latter assertion uses the hypotheses of
\cite[Theorem~3.20]{LiSuWu2026}, including the cuspidal global
eigensystem and the absolute irreducibility of the local residual
representation.
In either case, $\Tau=\JL(\rho)$.

Recently, Pengcheng Ding \cite{Ding2026} proves a strong finite length result 
about locally analytic vectors in a Banach representation.
In all cases considered here, $\Tau$ is an admissible Banach
representation of $D^\times$, and $\Tau^{\mathrm{an}}$ is an
admissible locally analytic representation of $D^\times$ with a
central character and, by \cite[Corollary~8.16]{DPS23}, an
infinitesimal character. Hence, by \cite[Theorem~1.1]{Ding2026},
$\Tau^{\mathrm{an}}$ has finite length.

In this article, we prove the irreducibility of $\Tau^{\mathrm{an}}$
when $\rho$ is irreducible and $\bar{\rho}$ satisfies the genericity
assumptions in \cite[Definition~4.4]{HuWangPadic}; see also
Hypothesis~\ref{hypothesis}.

\begin{mainresult}[Corollary~\ref{thm:tau-analytic-irreducible}]
Under Hypothesis~\ref{hypothesis}, the admissible locally analytic
representation $\Tau^{\mathrm{an}}$ of $D^\times$ is topologically
irreducible.
\end{mainresult}
To prove the result, 
we associate to each $D_n(G, E)$ module $M$ the residual symbol module $\mathsf{S}(M)$, defined in section~\ref{sec:distribution-symbols}.
when $G=U_D^1$, it is closely related
to the $\mathfrak{m}_D$-adic grading of a modulo $\varpi_E$ representation.

More generally, we prove the following non-existence theorem. Recall the ideal $J=(yz, zy)$:
\begin{mainresult}[Theorem~\ref{thm:no-two}]
For all sufficiently large $n$ satisfying
\[
 n\geq1+\frac{p+4}{2}+\frac1e,
\]
there is no finite $\mathscr D_{n,\zeta}$-module $M_n$, i.e. a $D_n(D^\times, E)$-module with central character $\zeta$, such that
$\mathsf S(M_n)$ is annihilated by a power of $J$ and
$\mu(M_n)=2$.
\end{mainresult}
And we apply it to the locally analytic vectors in admissible Banach representations.
\begin{mainresult}[Corollary~\ref{cor:banach-four-analytic-irreducible}]
      Let $\Pi$ be an admissible unitary Banach representation of
$D^\times$ with central character $\zeta^{-1}$.
Suppose that $\Pi'$ is pure over $E[\![U_D^1]\!]$ and
$\Pi^{\mathrm{an}}$ has an infinitesimal character. If
\[
 \Pi^{1\text{-}\mathrm{alg}}=0,\qquad \mu(\Pi)=4,
\]
then $\Pi^{\mathrm{an}}$ is topologically irreducible.
\end{mainresult}
Also, by \cite[Corollary~6.2]{LiuWang2025}, we may apply this criterion to these representations.

\paragraph{Structure of the paper}Here we give a brief outline of this article. 

Section~\ref{sec:two-filtrations} introduces filtered rings and
graded modules. We construct a notion of general residual symbol modules for certain filtered rings and modules,
and prove the well-definedness of characteristic cycles over a commutative quotient via the residual symbols.

In Section~\ref{sec:distributions}, we collect the distribution
algebras $D_n(G,E)$ and define their residual symbols
$\mathsf S(M_n)$. We then compare the residual symbols with the mod $p$ Iwasawa modules (Lemma~\ref{lem:large-radius-recovery}).
In Subsection~\ref{sec:analytic-to-residual}, we define ring homomorphisms
$D_n^\circ(G,E)\twoheadrightarrow k\llbracket G\rrbracket/J_d$,
where $D_n^\circ(G,E)$ is the unit ball. These maps give residual
representations of finite dimension from locally analytic representations (Proposition~\ref{prop:analytic-finite-jet-recovery}).

In Section~\ref{sec:equivariance}, we discuss the $D^\times$-action
and define the $J$-supported characteristic cycle. We then prove that
there is no admissible locally analytic representation whose 
residual symbol is annihilated by a power of $J$ and has
multiplicity 2 at every sufficiently large radius (Theorem~\ref{thm:no-two}). 

In Section~\ref{sec:application}, we first deduce an irreducibility
criterion for locally analytic vectors in certain Banach representations (Corollary~\ref{cor:banach-four-analytic-irreducible}).
We then apply the criterion to $\Tau(\rho)=\JL(\rho)/\JL(\rho)^{1-\mathrm{alg}}$, 
proving the main theorem in Corollary~\ref{thm:tau-analytic-irreducible}.

\paragraph{Acknowledgement} 
The author is deeply grateful to his advisor, Professor Yongquan Hu, for suggesting the problem studied in this paper and for his constant guidance and encouragement throughout this project. 
The author would also like to thank Benchao Su and Gabriel Dospinescu for helpful discussions. 
The author is convinced by all of them that their evidence showed the resulting $\Tau^{\mathrm{an}}$ should be irreducible.

\paragraph{The author's contribution and AI disclosure} The author used GPT-6 Astra for language editing and proofreading. GPT-6 Astra also proposed the construction of the map $\tau_d$ in Subsection~\ref{sec:analytic-to-residual} and provided a generalization of the finite-dimensional recovery lemma to infinite-dimensional representations, now stated as Proposition~\ref{prop:analytic-finite-jet-recovery}.
Before receiving this mathematical input from GPT-6 Astra, the author had developed the framework of residual symbols at all radii and, using the weaker finite-dimensional recovery result, proved that $\JL(\rho)^{\mathrm{an}}/\JL(\rho)^{1\text{-}\mathrm{alg}}$ has length at most two when the difference between the Hodge--Tate weights of $\rho$ is sufficiently large.
The author subsequently verified the construction and the generalized lemma, and showed that the stronger recovery result leads to the topological irreducibility theorem of this article.
The author takes full responsibility for the mathematical content and correctness of this article.

\section{General results on filtered rings and graded modules}
\label{sec:two-filtrations}
\subsection{Notation for filtered rings}
\label{sec:filtered-notions}
We first fix the notations used throughout this article, and mainly 
following \cite[\S\S1,3]{ST03} and
\cite[\S1.3 and \S2]{Sorgdrager2026}.
We fix a sufficiently large finite extension $E/\mathbb{Q}_p$, with the ring of integers $\cO_E$, uniformizer $\varpi_E$, and residue field $k$ as above 
and write $e=v_E(p)$ for the ramification index. 

We use decreasing exhaustive filtrations throughout, 
by which we mean that the filtration $F$ on a ring $A$ consists of additive
subgroups $(F^sA)_{s\in\mathbb R}$ satisfying
\[
 F^sA\supseteq F^tA\ (s\leq t),\qquad
 F^sA\,F^tA\subseteq F^{s+t}A,\qquad
 1\in F^0A,\qquad \bigcup_sF^sA=A.
\]
The filtration pieces $F^sA$ form a basis of open neighborhoods 
of zero. The filtration is \emph{separated} if $\bigcap_sF^sA=0$,
and \emph{complete and separated} if the canonical map of additive groups
\[
 A\longrightarrow\varprojlim_{s\to+\infty}A/F^sA
\]
is an isomorphism. We also put: 
\[
 F^{>s}A:=\bigcup_{t>s}F^tA,\qquad
 \gr_F^sA:=F^sA/F^{>s}A,\qquad
 \gr_FA:=\bigoplus_{s\in\mathbb R}\gr_F^sA.
\]

We use \emph{quasi-integral} filtrations, meaning that
$\gr_F^sA=0$ for $s\notin q^{-1}\mathbb Z$ for some integer $q\geq1$.
We also assume that there is a $\delta\in \mathbb{R}$ such that for all $s\in\mathbb{R}$
\[
 F^sA=F^{\delta\lceil s/\delta\rceil}A
 \quad(s\in\mathbb R).
\] We call such a filtration \emph{discrete}.
A homomorphism $f\colon A\to B$ of filtered rings is \emph{filtered} if 
$f(F^sA)\subset F^sB$ for every $s$.

We similarly define filtered modules over a filtered ring.
A filtered left $A$-module is an $A$-module $V$ with an exhaustive
decreasing filtration by additive subgroups such that
\[
 F^sA\,F^tV\subseteq F^{s+t}V\qquad(s,t\in\mathbb R).
\]
We use the preceding definitions of separation, completeness,
left-continuity, and discreteness also for $V$, and put
\[
 F^{>s}V:=\bigcup_{t>s}F^tV,\qquad
 \gr_FV:=\bigoplus_{s\in\mathbb R}F^sV/F^{>s}V.
\]
For $v\neq0$, set $\deg_F(v):=\sup\{s:v\in F^sV\}$.
When this supremum is finite and attained, the \emph{principal symbol}
of $v$ is defined as
\[
 \sigma_F(v):=v+F^{>\deg_F(v)}V
       \in\gr_F^{\deg_F(v)}V\setminus\{0\}.
\]
Set $\deg_F(0)=+\infty$ and $\sigma_F(0)=0$.
Note that the symbol $\sigma_F$ is neither additive nor multiplicative in general. 

There are several ways to induce filtrations.
For filtered modules $V_i$, a submodule $W\subseteq V$, and a graded
module $N$, we use the following conventions
unless stated otherwise.
\[
 \begin{aligned}
 &F^s\Bigl(\bigoplus_iV_i\Bigr):=\bigoplus_iF^sV_i,\quad
 F^sW:=W\cap F^sV,\quad
 &F^s(V/W):=(F^sV+W)/W,\\
 &F^s(V\langle d\rangle):=F^{s+d}V,
 &(N\langle d\rangle)_s:=N_{s+d}.
 \end{aligned}
\]
For a filtered homomorphism $A\to B$, the \emph{tensor filtration} is defined as 
\[
 F^s(B\otimes_AV):=
 \sum_{u+v\geq s}\im\bigl(F^uB\otimes_{F^0A}F^vV
                              \longrightarrow B\otimes_AV\bigr).
\]

An $A$-linear map $f:V\to V'$ is \emph{filtered} if
$f(F^sV)\subseteq F^sV'$ for every $s$. It is \emph{strict} if
\[
 f(F^sV)=f(V)\cap F^sV'\qquad(s\in\mathbb R).
\]
A short exact sequence of filtered modules is \emph{strict exact}
if both maps are strict.

We work with the good filtrations on modules in this article.
A filtration on $V$ is \emph{good} if there are finitely many elements
$v_i\in F^{\ell_i}V$ and real numbers $\ell_i$ such that
\[
 F^sV=\sum_{i=1}^qF^{s-\ell_i}A\,v_i
 \qquad(s\in\mathbb R).
\]

The main rings in this article have the properties introduced in
\cite[\S3]{ST03}. We use left Fr\'echet--Stein presentations
\[
 \mathscr A\xrightarrow{\sim}\varprojlim_n A_n
\]
in the sense of that reference. Each $A_n$ is a left noetherian
Banach $E$-algebra and is flat as a right $A_{n+1}$-module through
the transition map $A_{n+1}\to A_n$.

\label{def:coadmissible-stages}
A left $\mathscr A$-module $M$ is \emph{coadmissible} if there is an
inverse system of finite left $A_n$-modules $M_n$, whose transition
map $M_{n+1}\to M_n$ is $A_{n+1}$-linear via $A_{n+1}\to A_n$, such that
\[
 M\xrightarrow{\sim}\varprojlim_nM_n,\qquad
 A_n\otimes_{A_{n+1}}M_{n+1}\xrightarrow{\sim}M_n,
\]
where the second map is induced by the transition map and the first is
$\mathscr A$-linear. Equip each $M_n$ with its canonical quotient
Banach topology and $M$ with the projective-limit topology. We call
$M_n$ the $n$-th \emph{Banach stage} of this presentation.

\subsection{Residual symbols}
\label{sec:graded-modules}
Now we define the residual symbol modules. 
Let $A$ be a complete separated Zariskian filtered ring such that
\[
 \gr A=S[\epsilon^{\pm1}],
\]
where $S$ is a noetherian nonnegatively graded $k$-algebra and
$\epsilon$ is central of degree $\gamma>0$.
We assume that the ring filtration is quasi-integral, 
and so are the good filtrations on the modules below. 

Let $P$ be a finite graded $S[\epsilon^{\pm1}]$-module. A
\emph{graded integral lattice} in $P$ is a finite graded
$S[\epsilon]$-submodule $L\subset P$ such that
$L[\epsilon^{-1}]=P$. The $S[\epsilon]$-span of finitely many
homogeneous generators gives such a lattice. Multiplication by
$\epsilon$ is injective on $L$, and $L$ is bounded below in degree.

We put $\Fil_\epsilon^aP=\epsilon^aL$ for $a\in\mathbb Z$ and
write $\gr_\epsilon$ for the associated graded, also for
$S[\epsilon^{\pm1}]$ with its lattice $S[\epsilon]$.

\begin{definition}\label{def:symbol}
Let $M$ be a finite left $A$-module with a good filtration $F$,
and let $L\subset\gr_FM$ be a graded integral lattice.
Equip $\gr_FM$ with the $\epsilon$-filtration defined by $L$.
We define the corresponding residual symbol module as
\[
 \mathsf S(M)=\gr_\epsilon^0\gr_FM=L/\epsilon L.
\]
It is a finite graded $S$-module and depends on the filtration
and lattice.
\end{definition}

The next lemma gives compatibility of the residual symbol modules with subquotients.

\begin{lemma}
\label{lem:two-stage-strictness}
Let $M_2\subset M_1\subset M$ be finite $A$-modules, with a good
filtration $F$ on $M$. Give $M_1,M_2$ the induced filtrations
and $W=M_1/M_2$ the quotient filtration.
Choose an integral graded lattice $L\subset\gr_FM$.
On $\gr_FM_1,\gr_FM_2$ use the induced $\epsilon$-filtrations,
and on $\gr_FW$ use the quotient $\epsilon$-filtration. Thus
\[
 \begin{aligned}
 \Fil_\epsilon^0\gr_FM_i&=L\cap\gr_FM_i\quad(i=1,2),\\
 \Fil_\epsilon^0\gr_FW
 &=\im(L\cap\gr_FM_1\longrightarrow\gr_FW).
 \end{aligned}
\]
These are integral graded lattices. Both the $F$-filtrations and
the $\epsilon$-filtrations are good for the corresponding filtered rings.
The associated graded sequence is strict for the $\epsilon$-filtrations.
In particular,
\[
 0\longrightarrow\mathsf S(M_2)
 \longrightarrow\mathsf S(M_1)
 \longrightarrow\mathsf S(W)\longrightarrow0
\]
is exact. The induced lattices can be obtained from the principal symbols of suitable
filtered generators of the subquotient $W$.
\end{lemma}
\begin{proof}
Since the filtration on $A$ is Zariskian, the induced and quotient
$F$-filtrations are good. With these filtrations, we have exact sequence
\[
 0\longrightarrow\gr_FM_2\longrightarrow\gr_FM_1
 \longrightarrow\gr_FW\longrightarrow0.
\]
Since $S[\epsilon]$ is noetherian, the intersections
$L\cap\gr_FM_i$ are finite over $S[\epsilon]$. After inverting $\epsilon$ they give
$\gr_FM_i$. The image of the lattice $L$ in $\gr_FW$ is also finite and
spans $\gr_FW$ after inverting $\epsilon$, hence the image is a lattice and the $\epsilon$-filtration 
induced by it is a good filtration.

The maps are strict for the induced and quotient
$\epsilon$-filtrations by definition.
Thus the exact sequence of residual symbol 
follows by taking the degree $0$ part of the $\epsilon$-grading.

Finally, we choose finitely many homogeneous generators of $\mathsf{S}(W)$, say $\bar{m_1}, \dots, \bar{m_r}$, of degrees $\ell_i$, respectively. 
There are homogeneous elements $\tilde{m}_1, \dots, \tilde{m}_r\in\mathrm{Fil}_{\epsilon}^0\gr_FW$, again of degrees $\ell_i$ respectively, such that $\tilde{m}_i\ \mathrm{mod}\epsilon=\bar{m}_i$.
We can define a graded map of graded $S[\epsilon]$-modules \[\bigoplus_{i=1}^r S[\epsilon]\langle-\ell_i\rangle e_i\to\mathrm{Fil}_{\epsilon}^0\gr_FW, \qquad e_i\mapsto\tilde{m}_i\]
Thus by the graded Nakayama \cite[Lemma~6.2]{NegronPevtsova2023}, the map above is surjective as the reduction modulo $\epsilon$ is. Thus these elements $\tilde{m}_i$ generate $\mathrm{Fil}_{\epsilon}^0\gr_FW$ over $S[\epsilon]$, 
therefore these elements also generate $\gr_FW$ over $S[\epsilon^{\pm1}]$ by localization. 

In the same way, we can take $m_i\in W$ such that $\sigma_F(m_i)=\tilde{m}_i$, and they induce the map
\[\bigoplus_{i=1}^rA\langle-\ell_i\rangle e_i\to W, \qquad e_i\mapsto m_i.\]
Since $A$ is a Zariskian ring and the filtration $F$ on $W$ is complete, the map above is a strict surjection by \cite[Theorem~I.4.2.4.(5)]{LvO96}, as the associated graded map is surjective.
\end{proof}

\subsection{The characteristic cycle over a commutative quotient}
\label{sec:uniform}
We keep $A,S,\epsilon,\gamma$ and the filtrations of
Subsection~\ref{sec:graded-modules}. In this subsection we fix a homogeneous two-sided ideal
$J\subseteq S$ such that $S/J$ is commutative.

\paragraph{Residual symbols and their supports}
\label{sec:definitions}
We recall the support and cycle constructions in
\cite[Tags~02QW, 02QX and 02QZ]{Stacks} in the form we need.
For a finite graded $S$-module $N$ with $J^aN=0$, the module
\[
 \gr_JN=\bigoplus_{i=0}^{a-1}J^iN/J^{i+1}N
\]
is finite over $S/J$. We define its component cycle by
\begin{equation}\label{eq:component-cycle}
 \CC_J(N):=\sum_{\mathfrak p\in\Min\Supp_{S/J}(\gr_JN)}
 \length_{(S/J)_{\mathfrak p}}(\gr_JN)_{\mathfrak p}
 \, [V(\mathfrak p)].
\end{equation}
If $\dim\Supp_{S/J}(\gr_JN)\leq q$, we also define the $q$-dimensional cycle by
\begin{equation}\label{eq:q-cycle}
 \CC_q(N):=\sum_{\substack{\mathfrak p\in\Supp_{S/J}(\gr_JN)\\
                 \dim (S/J)/\mathfrak p=q}}
 \length_{(S/J)_{\mathfrak p}}(\gr_JN)_{\mathfrak p}
 \, [V(\mathfrak p)].
\end{equation}
Both sums have finitely many terms of finite length. We apply these
constructions to $N=\mathsf S(M)$ by the following finite d\'evissage;
compare \cite[Lemma~D.3.3]{HTT2008}.

\begin{lemma}\label{lem:local-length}
The category of all finite graded $S$-modules killed by a power of $J$ is a Serre
subcategory of all finite graded $S$-modules. For an exact sequence $0\to N'\to N\to N''\to0$ in this category, we have
\[
 \begin{aligned}
 \Supp_{S/J}(\gr_JN)
 &=\Supp_{S/J}(\gr_JN')\cup\Supp_{S/J}(\gr_JN''),\\
 \length_{(S/J)_{\mathfrak p}}(\gr_JN)_{\mathfrak p}
 &=\length_{(S/J)_{\mathfrak p}}(\gr_JN')_{\mathfrak p}
  +\length_{(S/J)_{\mathfrak p}}(\gr_JN'')_{\mathfrak p}
 \end{aligned}
\]
where the second equality holds whenever the lengths are finite.
These invariants can be computed from any finite filtration with factors killed by $J$ and
are unchanged by grading shifts. On modules with support dimension at most $q$,
$\CC_q$ is additive.
\end{lemma}
\begin{proof}
Submodules and quotients inherit $J$-power annihilation. For extension, we have
\[
 J^aN'=J^bN''=0\quad\Rightarrow\quad J^{a+b}N=0.
\]Thus the category is stable under subquotients and extensions. 

For two finite filtrations $F,G$ of $N$ with factors killed by $J$,
we use the modular law to identify the factors of common refinements:
\[
 \begin{aligned}
 \frac{F^{i+1}N+(F^iN\cap G^jN)}
      {F^{i+1}N+(F^iN\cap G^{j+1}N)}
 &\cong
 \frac{F^iN\cap G^jN}
      {(F^{i+1}N\cap G^jN)+(F^iN\cap G^{j+1}N)}\\
 &\cong
 \frac{G^{j+1}N+(G^jN\cap F^iN)}
      {G^{j+1}N+(G^jN\cap F^{i+1}N)}.
 \end{aligned}
\]
Refinement preserves the union of supports and the sum of localized
lengths. For the given exact sequence, we apply this comparison to
\[
 0\longrightarrow\frac{N'\cap J^iN}{N'\cap J^{i+1}N}
 \longrightarrow\frac{J^iN}{J^{i+1}N}
 \longrightarrow\frac{J^iN''}{J^{i+1}N''}
 \longrightarrow0.
\]
We localize and sum over $i$ to obtain the assertions.
\end{proof}

\paragraph{Independence of the integral graded lattice}
\label{sec:lattices}
Beyond graded modules over $S$, we need to compare different $S[\epsilon]$ lattices and their modulo $\epsilon$ reductions.
In order to do this, we apply the results of \cite[Theorem~1.1.17 and Corollary~1.1.18]{Ginzburg1998}
to $S[\epsilon]$ and its central regular element $\epsilon$.

\begin{proposition}\label{prop:lattice-invariance}
For a finite graded $S[\epsilon^{\pm1}]$-module $P$, the condition
$J^a(L/\epsilon L)=0$ for some $a$ is independent of the integral
graded lattice $L\subset P$. Under this condition, the support and the
localized lengths at minimal primes $\mathfrak p$,
\[
 \Supp_{S/J}\bigl(\gr_J(L/\epsilon L)\bigr),\qquad
 \length_{(S/J)_{\mathfrak p}}
 \bigl(\gr_J(L/\epsilon L)\bigr)_{\mathfrak p},
\]
and hence $\CC_J(L/\epsilon L)$, are independent of the choice of $L$.
\end{proposition}
\begin{proof}
We first consider the simple case $\epsilon L'\subseteq L\subseteq L'$ and denote by $T=L'/L$. Then we have
\[
 0\longrightarrow T\langle-\gamma\rangle\longrightarrow L/\epsilon L
 \longrightarrow L'/\epsilon L'\longrightarrow T\longrightarrow0.
\]
We apply Lemma~\ref{lem:local-length} to the two short exact sequences
through the image of the middle map. This proves the assertions in this case.

For arbitrary lattices, we can assume without loss of generality that $L\subseteq L'$, by replacing $L'$ by $L+L'$.
Choose $r\geq0$ with $\epsilon^rL'\subseteq L$ and put
\[
 \begin{gathered}
 K_i=L+\epsilon^iL'\quad(0\leq i\leq r),\\
 \epsilon K_i\subseteq K_{i+1}\subseteq K_i\quad(0\leq i<r).
 \end{gathered}
\]
We can then apply the preceding comparison successively from $K_0=L'$ to $K_r=L$.
\end{proof}

We use the same comparison for exact sequences.

\begin{proposition}\label{prop:graded-exact}
For an exact sequence $0\to P'\to P\to P''\to0$ of finite graded
$S[\epsilon^{\pm1}]$-modules and any integral graded lattices in them,
the middle reduction is annihilated by a power of $J$ if and only if
both end reductions are. Under this condition, their supports and
localized lengths satisfy the formulas in Lemma~\ref{lem:local-length}.
\end{proposition}
\begin{proof}
By Proposition~\ref{prop:lattice-invariance}, we may take
$L'=L\cap P'$ and $L''=\im(L\to P'')$. Since $\epsilon$ is regular on $L''$, the sequence
\[
 0\longrightarrow L'/\epsilon L'\longrightarrow L/\epsilon L
 \longrightarrow L''/\epsilon L''\longrightarrow0
\]
is exact. By applying Lemma~\ref{lem:local-length}, we get the desired assertion.
\end{proof}

\paragraph{Independence of the good module filtration}
\label{sec:filtrations}
Again beyond the choices of graded lattices, we need to compare different good filtrations on the original $M$. 
We recall the filtration comparison used in
\cite[\S3.3]{HuWangPadic} and \cite[proof of Lemma~3.81]{BHHMS2025};
see also \cite[Lemma~D.3.1]{HTT2008}. Recall that all the filtrations used are discrete with jumps
in $\delta\mathbb{Z}$.

\begin{lemma}\label{lem:filtration-step}
Let $F,G$ be good filtrations on a finite $A$-module $M$ such that
\[
 F^sM\subseteq G^sM\subseteq F^{s-\delta}M
 \qquad(s\in\delta\mathbb Z).
\]
Then $T=\bigoplus_sG^sM/F^sM$ is finite over $\gr A$ and there is an exact sequence
\[
 0\longrightarrow T\langle\delta\rangle
 \longrightarrow\gr_FM\longrightarrow\gr_GM
 \longrightarrow T\longrightarrow0.
\]
For the lattice reductions of $\gr_FM$ and $\gr_GM$, one is annihilated
by a power of $J$ if and only if the other is. In this case, their supports
and generic lengths agree.
\end{lemma}
\begin{proof}
The inclusion $\Fil^{a+\delta}A\,G^sM\subseteq F^{a+s}M$
defines the $\gr A$-action on $T$. In degree $s$, the natural map
$\gr_FM\to\gr_GM$ has kernel $T_{s+\delta}$ and cokernel $T_s$.
Thus the sequence is exact and $T$ is finite.

Apply Proposition~\ref{prop:graded-exact} to its two short exact
sequences and use invariance under grading shifts.
\end{proof}

We join arbitrary good filtrations by these one-step comparisons.

\begin{lemma}\label{lem:chain}
Any two good filtrations on a finite $A$-module are connected by
filtration shifts and pairs as in Lemma~\ref{lem:filtration-step}.
\end{lemma}
\begin{proof}
By \cite[Lemma~I.5.3]{LvO96}, choose integers $a,b\geq0$ with
$G^{s+b\delta}M\subseteq F^sM\subseteq G^{s-a\delta}M$ and put
\[
 H_j^sM=F^sM+G^{s+j\delta}M\qquad(-a\leq j\leq b).
\]
These filtrations are good, with
\[
 \begin{gathered}
 H_{-a}^sM=G^{s-a\delta}M,\qquad H_b^sM=F^sM,\\
 H_{j+1}^sM\subseteq H_j^sM\subseteq H_{j+1}^{s-\delta}M
 \quad(-a\leq j<b).
 \end{gathered}
\]
This is the required chain.
\end{proof}

\paragraph{Well-definedness of the characteristic cycle}
\label{sec:main}
We now combine the lattice and filtration comparisons.

\begin{theorem}\label{thm:main}
For a finite $A$-module $M$, the condition $J^a\mathsf S(M)=0$
for some $a$ is independent of the good filtration and integral
graded lattice. Under this condition, so are
\[
 \Supp_{S/J}\bigl(\gr_J\mathsf S(M)\bigr),\qquad
 \CC(M):=\CC_J(\mathsf S(M)).
\]
\end{theorem}
\begin{proof}
By Proposition~\ref{prop:lattice-invariance}, the annihilation condition,
support, and generic lengths are independent of the integral graded lattice.
We then apply Lemmas~\ref{lem:filtration-step} and \ref{lem:chain}
to obtain independence of the good filtration.
\end{proof}

Moreover, the notion of characteristic cycle is additive with respect to exact sequences.

\begin{theorem}\label{thm:additivity}
The category of all finite $A$-modules with $J^a\mathsf S(M)=0$ for some $a$ is a
Serre subcategory of all finite $A$-modules.
For an exact sequence $0\to M'\to M\to M''\to0$ in it,
\[
 \Supp_{S/J}(\gr_J\mathsf S(M))
 =\Supp_{S/J}(\gr_J\mathsf S(M'))\cup\Supp_{S/J}(\gr_J\mathsf S(M'')).
\]
For each $q\geq0$, its modules of support dimension at most $q$
form a Serre subcategory, on which
\[
 \CC_q(M)=\CC_q(\mathsf S(M)),\qquad
 \CC_q(M)=\CC_q(M')+\CC_q(M'').
\]
If the characteristic support is pure of dimension $q$, then
$\CC_q(M)=\CC(M)$.
\end{theorem}
\begin{proof}
Take a good filtration on $M$ and the induced filtrations on
$M',M''$. We apply Proposition~\ref{prop:graded-exact} to the exact sequence
\[
 0\longrightarrow\gr_FM'\longrightarrow\gr_FM
 \longrightarrow\gr_FM''\longrightarrow0.
\]
The proposition gives us the $J$-power condition and the support and length
formulas. These prove the assertions about Serre categories and the additivity of $\CC_q$.

By Theorem~\ref{thm:main} above, the result is independent of these choices.
The last assertion then 
follows from the definitions.
\end{proof}

The same argument gives us the case with one filtration.

\begin{corollary}\label{cor:one-stage-cycle}
If $\gr A=S$ instead, with the same remaining hypotheses, then
Theorems~\ref{thm:main} and \ref{thm:additivity} hold with
$\gr_FM$ in place of $\mathsf S(M)$ and no lattice choice.
\end{corollary}
\begin{proof}
Apply Lemma~\ref{lem:local-length} directly to the exact sequence
in Lemma~\ref{lem:filtration-step}, and use Lemma~\ref{lem:chain}.
For additivity, use induced filtrations on exact sequences.
\end{proof}

\subsection{Gabber filtrations}
\label{sec:gabber-filtrations}

We recall the result on filtered rings used to obtain pure graded modules.
We mainly follow \cite[Section~6]{HuWangPadic}.
For a finite module $V$ over an Auslander regular ring $A$, put
\[
 j_A(V)=\inf\{i:\operatorname{Ext}^i_A(V,A)\ne0\},
 \qquad j_A(0)=+\infty.
\]
A module $V\ne0$ is called pure of grade $j$ if
$j_A(V')=j$ for every nonzero submodule $V'\subseteq V$.

Let $A$ have a complete separated Zariskian filtration, and assume
that $A$ and $\gr A$ are Auslander regular. After rescaling the indices by a common factor, 
we assume the rings and modules are $\mathbb{Z}$-filtered. 
Under our decreasing convention, the Rees ring and Rees module are written as
\[
 \mathcal R(A)=\bigoplus_{i\in\mathbb Z}F^iA\,u^{-i},
 \qquad
 \mathcal R_F(V)=\bigoplus_{i\in\mathbb Z}F^iV\,u^{-i},
 \qquad \deg u=-1.
\]
Here $u$ is central and regular, and
\[
 \mathcal R(A)/u\mathcal R(A)\cong\gr A,
 \qquad
 \mathcal R_F(V)/u\mathcal R_F(V)\cong\gr_FV.
\]
Moreover, $\mathcal R(A)$ is Auslander regular by
\cite[Chapter~III, Theorem~3.1.7]{LvO96}.

For a pure finite $A$-module $V$, two good filtrations $F,G$ are called
\emph{tamely close} if
\[
 j_{\mathcal R(A)}\left(
 \frac{\mathcal R_{F+G}(V)}{\mathcal R_{F\cap G}(V)}
 \right)\geq j_A(V)+2,
\]
where $(F+G)^iV=F^iV+G^iV$ and
$(F\cap G)^iV=F^iV\cap G^iV$.
This is an equivalence relation on good filtrations.
We now recall the theorem defining the Gabber filtration.

\begin{theorem}
\label{thm:gabber-filtration}
Under these hypotheses, let $V$ be pure of grade $j$.
Each tame equivalence class contains a unique good filtration $F$, called the Gabber filtration,
for which $\gr_FV$ is pure of grade $j$ over $\gr A$.
For every $G$ in that class,
\[
 G^iV\subseteq F^iV\qquad(i\in\mathbb Z).
\]
The kernel of the resulting map
\[
 \gr_GV\longrightarrow\gr_FV
\]
is the largest submodule of $\gr_GV$ of grade greater than $j$.
\end{theorem}

This is \cite[Definition~6.4 and Theorem~6.5]{HuWangPadic},
which recalls \cite[Theorem~5.23]{BjorkEkstrom1990}.

\section{Distribution algebras and residual symbols}
\label{sec:distributions}
\subsection{The rings and filtrations associated to groups}
\label{sec:group-filtrations}
We begin with a compact valued $p$-adic group, its ring of continuous distributions (i.e. the Iwasawa algebra), and its ring of locally analytic distributions. 
We use the valuations and ordered expansions of \cite[\S4]{ST03}
and \cite[\S3]{Sorgdrager2026}; the identities for our radius indexing
are recorded below.

Let $G$ be a pro-$p$ analytic group of dimension $d$. 
As in \cite[III, \S2.1.2]{Lazard1965}, a function $\omega\colon G\to (\frac{1}{p-1}, +\infty]$ is called a $p$-adic valuation if it satisfies
\begin{itemize}
      \item $\omega(g)=+\infty$ iff $g=\mathrm{id}$, 
      \item $\omega(gh^{-1})\ge\min(\omega(g), \omega(h))$, 
      \item $\omega(ghg^{-1}h^{-1})\ge\omega(g)+\omega(h)$,
      \item $\omega(g^p)=\omega(g)+1$.
\end{itemize}
The valuation defines a sequence of subgroups $G_{\ge s}\coloneqq\{g\in G, \omega(g)\ge s\}$ and also $G_{> s}\coloneqq\{g\in G, \omega(g)> s\}$.
All of them are normal subgroups of $G$. 

By \cite[III, \S2.3.6]{Lazard1965}, $\omega$ also defines a filtration $\mathrm{Fil}_\omega$ on $k[G]$ by
\[\mathrm{Fil}_{\omega}^sk[G]=\mathrm{span}_k\{\prod_{i=1}^m(h_i-1): m\ge0, \sum\omega(h_i)\ge s\}.\]
There is an ordered basis $g_1,\dots, g_d\in G$ together with a homeomorphism
\[\begin{aligned}
\mathbb{Z}_p^d&\xrightarrow{\cong}G\\
(x_1,\dots,x_d)&\mapsto g_1^{x_1}\cdot\cdots\cdot g_d^{x_d}
\end{aligned}\]
such that $\omega(g_1^{x_1}\cdots g_d^{x_d})=\min(\omega(g_i)+v_p(x_i))$. We assume that the weights $\omega(g_i)$ are rational.
From now on we fix such a choice of ordered basis and denote $b_i=g_i-1$ and $\omega_i=\omega(g_i)$. We write $\underline{b}^{\underline{\alpha}}\coloneqq\prod_{i=1}^{d}(g_i-1)^{\alpha_i}$ for ordered product. 

The modulo $\varpi_E$ Iwasawa algebra $k\llbracket G\rrbracket\coloneqq\lim_{N<G}k[G/N]$ can be regarded as the completion of the filtered ring $(k[G], \mathrm{Fil}_\omega)$. 
Every element in $k\llbracket G\rrbracket$ has a unique expansion $\sum c_{\underline\alpha}\underline b^{\underline{\alpha}}$ for some sequence $\{c_{\underline\alpha}\in k\}$.
For a multi-index $\underline \alpha=(\alpha_1,\dots,\alpha_d)\in\mathbb{N}^d$, the weighted sum $\omega(\underline\alpha)=\sum\alpha_i\omega_i$ determines the filtration degree of the ordered monomials
\[    \deg_\omega(\underline{b}^{\underline{\alpha}})=\omega(\underline\alpha)=\sum\omega_i\alpha_i,\qquad
      \deg_\omega(\sum c_{\underline\alpha}\underline b^{\underline{\alpha}})=\min_{c_\alpha\neq0}(\sum\omega_i\alpha_i).\]

Define the graded $\mathbb F_p[\epsilon_p]$-Lie algebra
\[
 \gr_\omega G=\bigoplus_{s>0}G_{\geq s}/G_{>s},
\]
with Lie bracket and $\epsilon_p$ action
\[
 [\bar g,\bar h]=\overline{ghg^{-1}h^{-1}}
       \in G_{\geq s+t}/G_{>s+t},\qquad
 \epsilon_p\bar g=\overline{g^p}
       \in G_{\geq s+1}/G_{>s+1}
\]
for $\bar g\in\gr_\omega^sG$, $\bar h\in\gr_\omega^tG$.
The coefficient specialization in \cite[\S3.4, Remark~3.1]{Sorgdrager2026}
shows that
\[
 \gr_\omega^*k\llbracket G\rrbracket
 \cong k\otimes_{\mathbb F_p[\epsilon_p]}
       U_{\mathbb F_p[\epsilon_p]}(\gr_\omega G),
 \qquad \epsilon_p\longmapsto0.
\]
In particular, the Lie algebra enters this enveloping algebra via
$k\otimes_{\mathbb F_p[\epsilon_p]}\gr_\omega G$.

We now pass to integral Iwasawa algebra $\cO_E\llbracket G\rrbracket\coloneqq\lim_{N<G}\cO_E[G/N]$. We follow the same steps, and take $p$ into consideration. 
Since $(G, \omega)$ is a valued group, every elements in $\cO_E\llbracket G\rrbracket$ can be written as $\sum c_{\underline\alpha}\underline b^{\underline{\alpha}}$ for some sequence $\{c_{\underline\alpha}\in \cO_E\}$.
Now we define the $n$-th deformation degree and filtration $\mathrm{Fil}_n$ as 
\[\begin{aligned}
      \deg_n(\sum c_{\underline\alpha}\underline b^{\underline{\alpha}})&=\min_{c_{\underline\alpha}\neq0}(nv_E(c_{\underline\alpha})+\sum\omega_i\alpha_i)\\
      \mathrm{Fil}_n^r&\coloneqq\{f\in\cO_E\llbracket G\rrbracket, \deg_n(f)\ge r\}.
\end{aligned}\]
where $v_E$ is the valuation on $E$ normalized by $v_E(\varpi_E)=1$ and $n\in\mathbb{Q}_{>0}$.

Since $\cO_E\llbracket G\rrbracket$ is a complete local ring, all filtrations above are topologically equivalent to the maximal-ideal-adic filtration. Hence $\cO_E\llbracket G\rrbracket$ is complete under all $\mathrm{Fil}_n$.

What is interesting is the tensor product filtration on $E\llbracket G\rrbracket=\cO_E\llbracket G\rrbracket\otimes_{\cO_E}E$, and it can be reformulated as follows:
For each element $f$ in $E\llbracket G\rrbracket$, there is a unique expansion $f=\sum c_{\underline\alpha}\underline b^{\underline{\alpha}}$ for some bounded sequence $\{c_{\underline{\alpha}}\in E\}$, i.e. there is an $m\in\mathbb{Z}$ such that for all $\underline{\alpha}$, $c_{\underline{\alpha}}\in\varpi_E^{m}\cO_E$. 
The same formula extends $\deg_n$ and the filtration
$\mathrm{Fil}_n$ to $E\llbracket G\rrbracket$.

This filtration $\mathrm{Fil}_n$ is no longer complete when $d>0$, and the completion of $(E\llbracket G\rrbracket, \mathrm{Fil}_n)$ is exactly the distribution algebra $D_{p^{-{v_p(\varpi_E)}/{n}}}(G, E)$ with the filtration defined in \cite{ST03}. 
As a recall and check, we write down the precise elements in the completed ring:
\[ (E\llbracket G\rrbracket, \mathrm{Fil}_n)\widehat{  }=\{f=\sum c_{\underline\alpha}\underline b^{\underline{\alpha}}, nv_E(c_{\underline{\alpha}})+\omega(\underline{\alpha})\to\infty\text{ as } \omega(\underline\alpha)\to\infty\}\]

From now on we denote by $D_{n}$ or $D_n(G, E)$ the ring $(E\llbracket G\rrbracket, \mathrm{Fil}_n)\widehat{  }$ displayed above.
We also denote the extended filtration by $\mathrm{Fil}_n$; it is defined by the same function $\deg_n$. Thus the natural injections 
\[(\cO_E\llbracket G\rrbracket, \mathrm{Fil}_n)\to(E\llbracket G\rrbracket, \mathrm{Fil}_n)\to (D_n, \mathrm{Fil}_n)\]
are filtered, and the filtrations on their sources are induced from $(D_n, \mathrm{Fil}_n)$.

To simplify the notation, we ignore the fixed $G$ later and denote by \[\Lambda ^\circ=\cO_E\llbracket G\rrbracket,\quad \Lambda =E\llbracket G\rrbracket,\quad \overline{\Lambda} =k\llbracket G\rrbracket,\quad S \coloneqq\gr_\omega^*\overline{\Lambda} .\]

\begin{remark}
The norm of $f\in D_n(G,E)$ in \cite{ST03} is given as
\[
 p^{-\deg_n(f)/(ne)}.
\]
Thus $\deg_n$ gives the same filtration, up to reindexing.

      The reason we choose to index in this way is that we want to make the group operators in the same filtration indexes, and we can then construct in the next subsection about the doubly graded rings and modules, which will have the same grading base ring.
\end{remark}

\begin{lemma}\label{lem:graded-distribution-ring}
Assume $ne>1$, and put $\epsilon=\sigma_n(\varpi_E)$, of degree $n$.
There are isomorphisms of graded $k$-algebras
\[
 \gr_n^{*}\Lambda^\circ\cong S[\epsilon],\qquad
 \gr_n^{*}\Lambda\cong\gr_n^{*}D_n(G,E)
       \cong S[\epsilon^{\pm1}].
\]
\end{lemma}

\begin{proof}
Write
\[
 \underline b^{\underline\alpha}\underline b^{\underline\beta}
 =\sum_{\underline\gamma}c_{\underline\alpha\underline\beta,\underline\gamma}
                   \underline b^{\underline\gamma},
 \qquad c_{\underline\alpha\underline\beta,\underline\gamma}\in\mathbb Z_p.
\]
By \cite[Lemma~4.1]{ST03}, each nonzero coefficient $c=c_{\underline\alpha\underline\beta,\underline\gamma}$ satisfies
\[
 nv_E(c)+\omega(\underline\gamma)
 \geq\omega(\underline\alpha)+\omega(\underline\beta)
       +(ne-1)v_p(c).
\]
Since $ne>1$, terms with $p\mid c$ have strictly larger degree.
Thus the principal symbols of the ordered monomials multiply as in
$S$, and their bases give $\gr_n^{*}\Lambda^\circ\cong S[\epsilon]$.
Localizing at $\varpi_E$ and applying \cite[Lemma~4.3]{ST03}
gives the remaining isomorphisms.
\end{proof}

The rings $S$, $\overline\Lambda$,
$\Lambda$, and $D_n(G,E)$ for $ne>1$ are Auslander regular by
\cite[Chapter~III, Theorem~2.2.5]{LvO96}.
For $\Lambda$, we can also use the $\varpi_E$-adic filtration
\[
 F^i\Lambda=\varpi_E^i\Lambda^\circ,\qquad
 \gr_F\Lambda=\overline\Lambda[\epsilon^{\pm1}],\qquad
 \deg\epsilon=1.
\]
This filtration and the radius filtration on $D_n(G,E)$ are
complete and quasi-integral with noetherian associated graded
rings, hence Zariskian. Thus the results of
Subsection~\ref{sec:uniform} apply to $D_n(G,E)$.

Furthermore, if $(G,\omega)$ satisfies (HYP) of
\cite[\S4]{ST03}, then
\[
 S=k[x_1,\ldots,x_d],\qquad\deg x_i=\omega_i,
\]
by \cite[Lemma~4.4 and Theorem~4.5]{ST03}, and we may take $J=0$ in this case.

\subsection{The \texorpdfstring{$\epsilon$}{epsilon}--filtration and residual symbols}
\label{sec:distribution-symbols}
The Lemma~\ref{lem:graded-distribution-ring} above suggests us to separate the structure of $S$ out from a graded $S[\epsilon^{\pm1}]$ module. 
And we've prepared in Section~\ref{sec:graded-modules}.
We now apply these results to the situation
\[
 R=\gr_n^*D_n=S[\epsilon^{\pm1}],\qquad
 R^\circ=\gr_n^*\Lambda^\circ=S[\epsilon],\qquad \gamma=n.
\]
The $\epsilon$-filtration is a filtration of the associated graded
ring, distinct from the construction in
\cite[Sections~3.1--3.2]{ArdakovWadsley2013}. We use the following notations
\[
 \Fil_\epsilon^a\gr_n^*D_n=\epsilon^aS[\epsilon],
 \qquad \gr_\epsilon^0\gr_n^*D_n=S.
\]

Let $M_n$ be a finite $D_n$-module with a good filtration defined
by filtered generators $m_1,\ldots,m_s$. They determine the lattice
\[
 \Fil_\epsilon^0\gr_n^*M_n
 =\im\left(
 \bigoplus_iS[\epsilon]\langle-\deg(m_i)\rangle
 \longrightarrow\gr_n^*M_n
 \right).
\]
By Definition~\ref{def:symbol}, the corresponding residual symbol is given by
\[
 \mathsf S(M_n)=\gr_\epsilon^0\gr_n^*M_n.
\]

\begin{remark}
The lattice $\Fil_\epsilon^0\gr_n^*M_n$ depends on the filtered
presentation. Two sets of generators may give the same radius
filtration but different lattices and residual symbols.
We fix the presentation when comparing these symbols.
By Lemma~\ref{lem:two-stage-strictness}, suitable generators of
subquotients give the induced and quotient lattices.
\end{remark}

The following example compares different radii, and suggests that we can't have
a well-defined notion of characteristic cycle and multiplicity of a coadmissible module
across all (or close enough to 1) radii, at least not in the na\"ive sense. 

\begin{example}\label{ex:radius-dependent-support}
The support of the residual symbol may change with the radius.
In this example we take $J=0$ and $G=\mathbb Z_p^2$ with its standard valuation,
$E=\mathbb Q_p$, and $k=\mathbb F_p$. Let $g_1,g_2$ be an ordered
basis, put $b_i=g_i-1$, and write
\[
 S=\gr_\omega^*k\llbracket G\rrbracket=k[x,y],\qquad
 x=\sigma_\omega(b_1),\quad y=\sigma_\omega(b_2),\quad
 \deg x=\deg y=1.
\]
For each $m\geq1$, choose a homogeneous polynomial
$P_m(x,y)\in\mathbb Z_p[x,y]$ of degree $m^2$ whose reduction
$\overline P_m\in k[x,y]$ is irreducible of the same degree.
Such polynomials are obtained by lifting the coefficients of
homogenized monic irreducible polynomials over $\mathbb F_p$.
Set
\[
 f=\sum_{m\geq1}p^{-2m}P_m(b_1,b_2),\qquad
 M=D(G,\mathbb Q_p)/(f).
\]
For the radius indexing used in this article, we have
\[
 \mathsf S(M_n)\cong S/(\overline P_n),\qquad
 \Supp_S\mathsf S(M_n)=V(\overline P_n)
 \quad(n\geq2),
\]
where the residual symbol is taken with the quotient filtration
and its cyclic integral graded lattice. The supports are pairwise
distinct.
\end{example}
\begin{proof}
For every fixed radius index $N$, the coefficients of $P_m$ are
integral and at least one is a unit, so
\[
 \deg_N\bigl(p^{-2m}P_m(b_1,b_2)\bigr)=m^2-2Nm
 \longrightarrow+\infty\qquad(m\to\infty).
\]
Thus the series converges in every $D_N(G,\mathbb Q_p)$ and
defines $f\in D(G,\mathbb Q_p)$. The finitely presented module
$M$ is coadmissible by \cite[Corollary~3.4(v)]{ST03}, and
\[
 M_N=D_N(G,\mathbb Q_p)\otimes_{D(G,\mathbb Q_p)}M
     \cong D_N(G,\mathbb Q_p)/(f).
\]
At $N=n\geq2$, the equality
\[
 m^2-2nm=(m-n)^2-n^2
\]
shows that the minimum occurs uniquely at $m=n$. Consequently
\[
 \deg_n(f)=-n^2,\qquad
 \sigma_n(f)=\epsilon^{-2n}\overline P_n(x,y),\qquad
 \deg\epsilon=n.
\]
The radius norm is multiplicative, so multiplication by $f$ is
strict with its degree shift. With the quotient filtration,
\[
 \gr_n^*M_n
 \cong S[\epsilon^{\pm1}]/(\overline P_n).
\]
The image of $S[\epsilon]$ is the integral graded lattice
$(S/(\overline P_n))[\epsilon]$. Reducing it modulo $\epsilon$
gives us the stated residual symbol and its support. Since the irreducible
polynomials $\overline P_n$ have distinct degrees, their prime
ideals and supports are pairwise distinct. Moreover,
\[
 \CC(M_n)=[V(\overline P_n)].
\]
Thus neither the support nor the characteristic cycle need
stabilize as the radius index grows. 
\end{proof}

\subsection{Large-radius recovery for finite Iwasawa modules}
\label{sec:recovery}

Let $M^\circ$ be a finitely generated, $\varpi_E$-torsion-free left
$\Lambda ^\circ=\cO_E\llbracket G\rrbracket$-module and $M=M^\circ[1/\varpi_E]$.  Define, on separate
levels,
\[
 M_{\infty}
 :=D(G,E)\otimes_{\cO_E\llbracket G\rrbracket}M^\circ,\qquad
 M_n
 :=D_n(G, E)\otimes_{\cO_E\llbracket G\rrbracket}M^\circ.
\]
When $D(G,E)=\varprojlim_nD_n$ is a Fr\'echet--Stein presentation,
these modules satisfy, by \cite[Theorem~A and Corollaries~3.1,~3.4(v)]{ST03},
\[
 M_{\infty}\xrightarrow{\sim}\varprojlim_nM_n,\qquad
 D_n\otimes_{D(G,E)}M_{\infty}\xrightarrow{\sim}M_n.
\]
Choose generators $m_1,\ldots,m_s\in M^\circ$ and fixed shifts
$\ell_1,\ldots,\ell_s$.  Give $M^\circ$ the quotient radius filtration $\mathrm{Fil}_{Q,n}^*$ induced
by
\[
 (F^\circ,\mathrm{Fil}_{\underline{\ell}, n}^*)=(\bigoplus_{i=1}^s
 \cO_E\llbracket G\rrbracket\langle-\ell_i\rangle, \mathrm{Fil}_n^*)
 \twoheadrightarrow (M^\circ, \mathrm{Fil}_{Q, n}^*),
\]
where $\mathrm{Fil}_{n}^*$ is the subspace filtration from $\cO_E\llbracket G\rrbracket\subset D_n(G, E)$,
and give $M_n$ the quotient filtration induced by the
base-changed presentation
\[
 \bigoplus_{i=1}^sD_n\langle-\ell_i\rangle
 \twoheadrightarrow M_n.
\]

\begin{lemma}
\label{lem:large-radius-recovery}
With the preceding notation:
\begin{enumerate}
 \item for every $n$ with $ne>1$ at which the displayed presentations
       define good filtrations, there is a natural graded
       exact sequence
       \[0\to\gr_{Q,n}^{*}M^\circ[\epsilon^{\infty}]\to
        \gr_{Q,n}^{*}M^\circ\twoheadrightarrow
        \Fil_\epsilon^0\gr_n^{*}M_n\to0,
       \]where the first term means the $\epsilon$-power-torsion in $\gr_{Q,n}^{*}M^\circ$;
 \item for $n\gg0$, $\gr_{Q,n}^{*}M^\circ$ is
       $\epsilon$-torsion-free and in particular the map above is an isomorphism:
       \[
        \gr_{Q,n}^{*}M^\circ
        \xrightarrow{\sim}
        \Fil_\epsilon^0\gr_n^{*}M_n;
       \]
 \item for the same $n$, multiplication by $\varpi_E$ is strict and
       \[
        \gr_{Q,n}^{*}(M^\circ/\varpi_EM^\circ)
        \xrightarrow{\sim}
        \gr_\epsilon^0\gr_n^{*}M_n,
       \]
       and the left side is independent of $n$ when $n\gg0$.
\end{enumerate}
\end{lemma}

\begin{proof}
Put $R^\circ=\ker(F^\circ\to M^\circ)$ and
$F=F^\circ[1/\varpi_E]$, and give $M$ the quotient
filtration from $F$. By \cite[Lemmas~2.1--2.3]{Sorgdrager2026},
localization and completion give the commutative diagram with
surjective rows
\[
\begin{tikzcd}
 \displaystyle\bigoplus_iS[\epsilon]\langle-\ell_i\rangle
   \arrow[r,two heads]\arrow[d]
 &\gr_{Q,n}^*M^\circ\arrow[r]\arrow[d,"\lambda_n"]&0\\
 \displaystyle\bigoplus_iS[\epsilon^{\pm1}]\langle-\ell_i\rangle
   \arrow[r,two heads]\arrow[d,"\sim"']
 &\gr_{Q,n}^*M\arrow[r]\arrow[d,"\sim"]&0\\
 \displaystyle\bigoplus_i\gr_n^*D_n\langle-\ell_i\rangle
   \arrow[r,two heads]
 &\gr_n^*M_n\arrow[r]&0.
\end{tikzcd}
\]
Here $\lambda_n$ is localization at $\epsilon$, followed by
$(\gr_{Q,n}^*M^\circ)[\epsilon^{-1}]
 \xrightarrow{\sim}\gr_{Q,n}^*M$.
Thus the composite $\iota_n:\gr_{Q,n}^*M^\circ\to\gr_n^*M_n$ satisfies
\[
 \ker(\iota_n)=\gr_{Q,n}^*M^\circ[\epsilon^\infty],\qquad
 \im(\iota_n)=\Fil_\epsilon^0\gr_n^*M_n,
\]
which proves (1).

For (2), apply the argument of
\cite[Proposition~5.1]{Sorgdrager2026}. It also applies to our fixed
shifts: since $M^\circ$ is $\varpi_E$-torsion-free,
\[
 \overline R:=R^\circ/\varpi_ER^\circ
 \hookrightarrow F^\circ/\varpi_EF^\circ
       =\bigoplus_i\overline\Lambda\langle-\ell_i\rangle.
\]
Give $\overline R$ the subspace filtration in this display. Choose
lifts $v_j\in R^\circ$ of elements whose principal symbols generate
$\gr_\omega^*\overline R$, and let $d_j$ be their degrees.
The radius estimate in that proof becomes
\[
 n+\min_i\ell_i>\max_jd_j
 \quad\Longrightarrow\quad \deg_n(v_j)=d_j,
\]
because each coefficient divisible by $\varpi_E$ has degree at least
$n+\min_i\ell_i$. Hence $\gr_{Q,n}^*M^\circ$ is
$\epsilon$-torsion-free for all sufficiently large $n$, and (1) gives (2).

For such $n$, multiplication by $\varpi_E$ is strict. Taking associated
gradeds of
\[
 0\longrightarrow M^\circ\langle-n\rangle
 \xrightarrow{\varpi_E}M^\circ
 \longrightarrow M^\circ/\varpi_EM^\circ\longrightarrow0
\]
and using (2) gives the commutative diagram with exact rows
\[
\begin{tikzcd}[column sep=small,row sep=large]
 0\arrow[r]
 &\gr_{Q,n}^*M^\circ\langle-n\rangle
   \arrow[r,"\epsilon"]\arrow[d,"\sim"']
 &\gr_{Q,n}^*M^\circ\arrow[r]\arrow[d,"\sim"]
 &\gr_{Q,n}^*(M^\circ/\varpi_EM^\circ)
   \arrow[r]\arrow[d,dashed,"\sim"]&0\\
 0\arrow[r]
 &\Fil_\epsilon^0\gr_n^*M_n\langle-n\rangle
   \arrow[r,"\epsilon"]
 &\Fil_\epsilon^0\gr_n^*M_n\arrow[r]
 &\mathsf S(M_n)\arrow[r]&0.
\end{tikzcd}
\]
The quotient filtration on $M^\circ/\varpi_EM^\circ$ is induced from
$\bigoplus_i\overline\Lambda\langle-\ell_i\rangle$ with its fixed
$\omega$-filtration, so it is independent of $n$ in our indexing.
This proves (3), as in \cite[Proposition~2.1 and Lemma~2.4]{Sorgdrager2026}.
\end{proof}

\subsection{From analytic representations to residual representations}
\label{sec:analytic-to-residual}

We fix $n$ with $ne>1$ and a finite $D_n$-module $M_n$ with a complete
separated good filtration $\mathrm{Fil}_n^*M_n$ and an integral graded lattice
as in Subsection~\ref{sec:distribution-symbols}. We use its residual
symbol $\mathsf S(M_n)$ and put
\[
 \mathsf S(M_n)_{\geq d}
       :=\bigoplus_{s\geq d}\mathsf S(M_n)_s.
\]
We impose the following bound on the relation degrees.

\begin{hypothesis}
After a common grading shift, $\mathsf S(M_n)$ admits a graded presentation
\begin{equation}\label{eq:jet-bounded-presentation}
 \begin{gathered}
 0\longrightarrow I\longrightarrow S^{q}\longrightarrow\mathsf S(M_n)
     \longrightarrow0,\\
 I=\sum_iS\bar a_i,\qquad
 s_i:=\deg\bar a_i,\qquad 0\leq s_i\leq d'<n,
 \end{gathered}
\end{equation}
where $d'\geq0$, the free generators have degree zero, and the $\bar a_i$
are homogeneous.
\end{hypothesis}

\paragraph{The truncation map}
We first define the finite quotient used in the construction.
We fix an integer $d>0$ and denote
\[
 J_d:=\mathrm{Fil}_\omega^d\overline\Lambda,\qquad
 D_n^\circ:=\mathrm{Fil}_n^0D_n,\qquad M_n^\circ:=\mathrm{Fil}_n^0M_n.
\]
Then $\overline{\Lambda}/J_d$ has the finite $k$-basis
$\{\bar b^\alpha:\omega(\alpha)<d\}$. For $n>d'+d$, define
\[
 \tau_d:\mathrm{Fil}_n^{-d'}D_n\longrightarrow \overline{\Lambda}/J_d,
 \qquad
 \sum_\alpha c_\alpha b^\alpha\longmapsto
 \sum_{\omega(\alpha)<d}\overline{c}_\alpha\,\bar b^\alpha.
\]
The map is well defined since
\[
 nv_E(c_\alpha)+\omega(\alpha)\geq-d',\quad\omega(\alpha)<d<n-d'
 \quad\Longrightarrow\quad v_E(c_\alpha)\geq0.
\]
The radius index $n$ is understood from the domain of $\tau_d$.

Put $\mathfrak a_{n,d}:=\ker(\tau_d|_{D_n^\circ})$. We will show below
that it is a two-sided ideal of $D_n^\circ$. The main result to be proved in this
subsection is the following.

\begin{proposition}
\label{prop:analytic-finite-jet-recovery}
Assume \eqref{eq:jet-bounded-presentation}. For $d>0$ and
$n\geq d'+d+1/e$, the finite $\overline\Lambda$-module
\[
 \overline{M}_{n,d}^{\circ}:=M_n^\circ/\mathfrak a_{n,d}M_n^\circ
 \cong(\overline\Lambda/J_d)\otimes_{D_n^\circ,\tau_d}M_n^\circ
\]
with its quotient filtration from $M_n^\circ$ satisfies
\[
 \gr_{Q,n}^*\overline{M}_{n,d}^{\circ}
 \cong\mathsf S(M_n)/\mathsf S(M_n)_{\geq d}.
\]
Here $\mathrm{Fil}_{Q,n}^s\overline{M}_{n,d}^{\circ}=\overline{M}_{n,d}^{\circ}$ for $s<0$, and
$\mathrm{Fil}_{Q,n}^s\overline{M}_{n,d}^{\circ}=\im(\mathrm{Fil}_n^sM_n\to \overline{M}_{n,d}^{\circ})$
for $s\geq0$.
\end{proposition}

Before proving the proposition, we first establish the coefficient
estimate for $\tau_d$.

\begin{lemma}\label{lem:analytic-jet-truncation}
Assume $n\geq d'+d+1/e$. Then $\tau_d|_{D_n^\circ}$ is a ring homomorphism and
\[
 \tau_d:\mathrm{Fil}_n^{-d'}D_n\twoheadrightarrow\overline\Lambda/J_d
\]
is a morphism of $D_n^\circ$-bimodules.
Moreover, $\tau_d(\mathrm{Fil}_n^dD_n)=0$, and
$\mathfrak a_{n,d}$ is a two-sided ideal of $D_n^\circ$.
\end{lemma}

\begin{proof}
By definition, $\tau_d$ is $\cO_E$-linear and
$\tau_d(\mathrm{Fil}_n^dD_n)=0$. For
$f\in\mathrm{Fil}_n^{-d'}D_n$ and $g\in D_n^\circ$, choose finite
monomial sums $f_0,g_0$ such that
\[
 f-f_0\in\mathrm{Fil}_n^dD_n,\qquad
 g-g_0\in\mathrm{Fil}_n^{d+d'}D_n
 \quad\Rightarrow\quad
 fg-f_0g_0\in\mathrm{Fil}_n^dD_n.
\]
By bilinearity, it suffices to take $f=\lambda b^\alpha$ and
$g=\mu b^\beta$ with $\lambda,\mu\ne0$ and to show
$\tau_d(f)\tau_d(g)=\tau_d(fg)$. Write
\[
 b^\alpha b^\beta=\sum_\gamma c_{\alpha\beta\gamma}b^\gamma.
\]
For a nonzero coefficient $c_{\alpha\beta\gamma}$, \cite[Lemma~4.1]{ST03} and the filtration give us the first estimation
\[
 \begin{gathered}
 v_E(c_{\alpha\beta\gamma})\geq0,\qquad
 \omega(\gamma)\geq\omega(\alpha)+\omega(\beta)
                   -\frac1e v_E(c_{\alpha\beta\gamma}),\\
 n\bigl(v_E(\lambda)+v_E(\mu)+v_E(c_{\alpha\beta\gamma})\bigr)
       +\omega(\gamma)\geq-d'.
 \end{gathered}
\]
If $\omega(\gamma)<d<n-d'$, the last inequality implies
\[
 v_E(\lambda)+v_E(\mu)+v_E(c_{\alpha\beta\gamma})\geq0.
\]
For a term surviving reduction, this sum is zero. Hence
\[
 \begin{aligned}
 v_E(c_{\alpha\beta\gamma})
 &=-v_E(\lambda)-v_E(\mu)\\
 &\leq\max\{-v_E(\lambda),0\}+\max\{-v_E(\mu),0\},\\
 \omega(\gamma)
 &\geq n\bigl(\max\{-v_E(\lambda),0\}+\max\{-v_E(\mu),0\}\bigr)
       -d'-\frac1e v_E(c_{\alpha\beta\gamma})\\
 &\geq(n-1/e)\bigl(\max\{-v_E(\lambda),0\}
                         +\max\{-v_E(\mu),0\}\bigr)-d'.
 \end{aligned}
\]
Consequently,
\[
 v_E(\lambda)<0\text{ or }v_E(\mu)<0
 \quad\Rightarrow\quad
 \omega(\gamma)\geq n-1/e-d'\geq d.
\]
In this case $\tau_d(f)\tau_d(g)=0=\tau_d(fg)$.
If $v_E(\lambda),v_E(\mu)\geq0$, both factors lie in $\Lambda^\circ$, and $\tau_d|_{\Lambda^\circ}$
is the usual reduction modulo $(\varpi_E, J_d)$. Reversing the factors gives us
\[
 \tau_d(fg)=\tau_d(f)\tau_d(g),\qquad
 \tau_d(gf)=\tau_d(g)\tau_d(f).
\]
The restriction to $D_n^\circ$ is therefore a unital ring map, and
surjectivity follows from $\Lambda^\circ\subset D_n^\circ$.

\end{proof}

We next lift the relations in \eqref{eq:jet-bounded-presentation}.
\begin{lemma}\label{lem:analytic-jet-relations}
There are graded isomorphisms
\[
 \Fil_\epsilon^0\gr_n^*M_n\cong\mathsf S(M_n)[\epsilon],\qquad
 \gr_n^*M_n\cong\mathsf S(M_n)[\epsilon^{\pm1}],
\]
and a strict filtered presentation whose free generators have degree zero:
\[
 0\longrightarrow K_n\longrightarrow F_n:=D_n^{q}
   \xrightarrow{\varphi_n}M_n\longrightarrow0,\qquad
 \gr_n^*K_n=I[\epsilon^{\pm1}].
\]
Each $\bar a_i$ lifts to $a_i\in K_n$ with
$\deg_n(a_i)=s_i$ and $\sigma_n(a_i)=\bar a_i$.
\end{lemma}

\begin{proof}
Write $\pi_n:S^{q}\twoheadrightarrow\mathsf S(M_n)$ for the
presentation in \eqref{eq:jet-bounded-presentation}, so $\ker\pi_n=I$.
Choose homogeneous degree-zero lifts of its generators in
$\Fil_\epsilon^0\gr_n^*M_n$. Graded Nakayama
\cite[Lemma~6.2]{NegronPevtsova2023}, applied to the cokernel,
gives us the surjective top map in the commutative diagram
\[
\begin{tikzcd}[column sep=large,row sep=large]
 S[\epsilon]^{q}\arrow[r,two heads,"\Phi_n"]
   \arrow[d,two heads,"\epsilon\mapsto0"']
 &\Fil_\epsilon^0\gr_n^*M_n
   \arrow[d,two heads,"\bmod\epsilon"]\\
 S^{q}\arrow[r,two heads,"\pi_n"']
 &\mathsf S(M_n).
\end{tikzcd}
\]
Its source is concentrated in nonnegative degrees, hence
$(\Fil_\epsilon^0\gr_n^*M_n)_{<0}=0$.
Since $\bar a_i\in I=\ker\pi_n$, and by viewing 
$\bar a_i$ as an element in $S[\epsilon]^q$, commutativity gives
\(
 \Phi_n(\bar a_i)\bmod\epsilon
   =\pi_n(\bar a_i)=0.
\)
Thus its image is divisible by $\epsilon$. Since $\Phi_n$ is a graded morphism 
of graded $S[\epsilon]$ modules and $s_i\leq d'<n$, more precisely we have
\[
 \Phi_n(\bar a_i)
 \in\epsilon(\Fil_\epsilon^0\gr_n^*M_n)_{s_i-n}=0.
\]
The generators of $I$ therefore map to zero, so
\[
 I[\epsilon]\subseteq\ker\Phi_n,\qquad
 \frac{S[\epsilon]^{q}}{I[\epsilon]}
       \cong\mathsf S(M_n)[\epsilon].
\]
Consequently $\Phi_n$ induces a surjection of $S[\epsilon]$-modules
\[
 \Psi_n:\mathsf S(M_n)[\epsilon]
       \twoheadrightarrow\Fil_\epsilon^0\gr_n^*M_n
\]
whose reduction modulo $\epsilon$ is the identity of $\mathsf S(M_n)$.
The target has no $\epsilon$-torsion, so reduction of the kernel
sequence gives us
\[
 (\ker\Psi_n)/\epsilon(\ker\Psi_n)=0
 \quad\Rightarrow\quad\ker\Psi_n=0
\]
by graded Nakayama. Thus also $\mathrm{ker}\ \Phi_n=I[\epsilon]$.

By inverting $\epsilon$, we get the 
graded isomorphism $\gr_nM_n\cong\mathsf{S}(M_n)[\epsilon^{\pm1}]$.
By complete filtered lifting to $M_n$, \cite[Theorem~I.4.2.4.(5)]{LvO96},
we get the strict presentation and the required lifts $a_i\in K_n$, as the last step in lemma~\ref{lem:two-stage-strictness}.
\end{proof}

\paragraph{Recovery of finite-dimensional representations and gradings}

Let us recollect the settings and notations for the rest of this subsection.
Fix $d>0$ and $n\geq d'+d+1/e$,
and use the presentation and chosen lifts of
Lemma~\ref{lem:analytic-jet-relations}. Put
\[
 F_n^\circ:=(D_n^\circ)^{q}=\bigoplus_iD_n^\circ\mathbf e_i,
 \qquad K_n^\circ:=K_n\cap F_n^\circ.
\]
Recall that for the kernel in the presentation,
\[
 \begin{gathered}
 \gr_n^*K_n=I[\epsilon^{\pm1}],\qquad
 I=(\gr_n^*K_n^\circ)\cap S^{q}=\sum_iS\bar a_i,\\
 \text{lifts }a_i\in K_n^\circ,\qquad
 \sigma_n(a_i)=\bar a_i,\qquad
 \deg_n(a_i)=s_i\in[0,d'].
 \end{gathered}
\]
Here $K_n^\circ$ has the induced filtration, and the intersection is taken
inside $S[\epsilon^{\pm1}]^{q}$. On $F_n^\circ$, the notation $\tau_d$
denotes the coefficientwise map
\[
 F_n^\circ\longrightarrow
 (\overline\Lambda/J_d)\otimes_{D_n^\circ,\tau_d}F_n^\circ
 \cong(\overline\Lambda/J_d)^{q},\qquad
 x\longmapsto1\otimes x.
\]
Its restriction to $K_n^\circ$ and the quotient map
$M_n^\circ\twoheadrightarrow\overline M_{n,d}^\circ$ are also denoted by
$\tau_d$. Since
\[
 (\varphi_n|_{F_n^\circ})^{-1}(\mathfrak a_{n,d}M_n^\circ)
 =K_n^\circ+\mathfrak a_{n,d}F_n^\circ,
\]
the kernel in the last row below is $\tau_d(K_n^\circ)$.
The rational presentation and its integral and residual forms
give the following commutative diagram with exact rows:
\[
\begin{tikzcd}[column sep=large,row sep=large]
 0\arrow[r]
 &K_n\arrow[r,hook]
 &F_n=D_n^{q}\arrow[r,two heads,"\varphi_n"]
 &M_n\arrow[r]&0\\
 0\arrow[r]
 &K_n^\circ\arrow[r,hook]\arrow[u,hook]
    \arrow[d,two heads,"\tau_d"']
 &F_n^\circ\arrow[r,two heads,"\varphi_n"]\arrow[u,hook]
    \arrow[d,two heads,"\tau_d"']
 &M_n^\circ\arrow[r]\arrow[u,hook]
    \arrow[d,two heads,"\tau_d"]&0\\
 0\arrow[r]
 &\tau_d(K_n^\circ)\arrow[r,hook]
 &(\overline\Lambda/J_d)^{q}\arrow[r,two heads]
 &\overline M_{n,d}^\circ\arrow[r]&0.
\end{tikzcd}
\]
In particular,
\[
 \overline{M}_{n,d}^{\circ}
 \cong(\overline\Lambda/J_d)^{q}/\tau_d(K_n^\circ).
\]
For $s\geq0$, coefficientwise truncation and the degree-zero free generators give us the following identities
\[
 \tau_d(\mathrm{Fil}_n^sD_n)
   =\mathrm{Fil}_\omega^s(\overline\Lambda/J_d),\qquad
 \mathrm{Fil}_{Q,n}^s\overline{M}_{n,d}^{\circ}
   =(\mathrm{Fil}_\omega^s\overline\Lambda)\overline{M}_{n,d}^{\circ}.
\]
We equip $\tau_d(K_n^\circ)$ with the induced subspace filtration
\[
 \mathrm{Fil}_\omega^s\tau_d(K_n^\circ)
 =\tau_d(K_n^\circ)\cap
   \mathrm{Fil}_\omega^s(\overline\Lambda/J_d)^{q}.
\]
For $0\leq j<d$, the inclusions of the integral modules induce
\[
 \gr_n^jK_n^\circ=\gr_n^jK_n,
 \qquad \gr_n^jF_n^\circ=\gr_n^jF_n,
 \qquad \gr_n^jM_n^\circ=\gr_n^jM_n.
\]
Using these identifications and
Lemma~\ref{lem:analytic-jet-relations}, we can write down the graded truncation explicitly as
\[
 \gr^j(\tau_d):
 \bigl(S[\epsilon^{\pm1}]^{q}\bigr)_j
 \twoheadrightarrow S_j^{q},\qquad
 \sum_h\epsilon^hv_h\longmapsto v_0,
 \qquad v_h\in S_{j-nh}^{q}.
\]
Here $S_j^{q}=\gr_\omega^j((\overline\Lambda/J_d)^{q})$.
Next, we prove surjectivity on the filtered relation modules.

\begin{lemma}\label{lem:analytic-jet-degree-lift} The map $\tau_d\colon K_n^\circ\to\tau_d(K_n^\circ)$ is strict with respect to the filtration defined above. 
That is, for $0\leq j<d$, we have a surjection
\[
 \tau_d:\Fil_n^jK_n^\circ
 \twoheadrightarrow\Fil_\omega^j\tau_d(K_n^\circ)=\tau_d(K_n^\circ)\cap\mathrm{Fil}_{\omega}^j(\overline{\Lambda}/J_d)^{q}.
\]
Consequently, we have a graded surjection for $0\le j<d$
\[
 \gr^j(\tau_d):\gr_n^jK_n=(I[\epsilon^{\pm1}])_j
 \twoheadrightarrow\gr_\omega^j\tau_d(K_n^\circ).
\]
\end{lemma}

\begin{proof}
We fix $0\leq j<d$ and
$0\neq\bar x\in\Fil_\omega^j\tau_d(K_n^\circ)$, and put
\[
 t:=\deg_\omega(\bar x),\qquad j\leq t<d.
\]
We choose a lift $x\in K_n^\circ$ of maximal degree $a=\deg_n(x)$.
Such a lift exists because the degrees are discrete and $0\leq a\leq t$.
Suppose $a<t$. Lemma~\ref{lem:analytic-jet-relations} gives us the expansion
\[
 \sigma_n(x)=\sum_{h\in\mathbb Z}\epsilon^h v_h,
 \qquad v_h\in I_{a-nh}.
\]
Since $I_{<0}=0$, the coefficient formula for graded truncation
gives us
\[
 \begin{aligned}
 h>0&\Rightarrow a-nh<0\Rightarrow v_h=0,\\
 a<t&\Rightarrow v_0=\gr^a(\tau_d)(\sigma_n(x))=0.
 \end{aligned}
\]
For $h<0$, we write
\[
 v_h=\sum_i f_{i,h}\bar a_i,
 \qquad f_{i,h}\in S_{a-nh-s_i}.
\]
We take $f_{i,h}=0$ when $v_h=0$, and put
\[
 \bar q_i:=\sum_{h<0}\epsilon^h f_{i,h},\qquad
 \sigma_n(x)=\sum_i\bar q_i\bar a_i.
\]
We can choose ordered termwise monomial lifts of $\bar{q}_i$ satisfying
$q_i\in\Fil_n^{a-s_i}D_n\subseteq\Fil_n^{-d'}D_n$.
Every monomial in $q_i$ has group weight at least
\[
 n+a-s_i\geq n-d'>d\quad\Rightarrow\quad\tau_d(q_i)=0.
\]
Lemma~\ref{lem:analytic-jet-truncation} therefore gives us
\[
 x':=x-\sum_iq_i a_i\in K_n^\circ,\qquad
 \tau_d(x')=\bar x,\qquad \deg_n(x')>a.
\]
This contradicts the choice of $x$ with maximal degree $a$. Hence
\[
 a=t\geq j\quad\Rightarrow\quad
 x\in\Fil_n^jK_n^\circ,\qquad\tau_d(x)=\bar x.
\]
The zero element lifts to zero, so we have the filtered surjection.

Finally, we use $\Fil_n^jK_n^\circ=\Fil_n^jK_n$ for $j\ge0$ in the
commutative diagram
\[
\begin{tikzcd}[column sep=large,row sep=large]
 \Fil_n^jK_n^\circ\arrow[r,two heads,"\tau_d"]\arrow[d,two heads]
 &\Fil_\omega^j\tau_d(K_n^\circ)\arrow[d,two heads]\\
 \gr_n^jK_n=(I[\epsilon^{\pm1}])_j
   \arrow[r,"\gr^j(\tau_d)"']
 &\gr_\omega^j\tau_d(K_n^\circ).
\end{tikzcd}
\]
The top and right maps are surjective, which gives us the graded surjection.
\end{proof}

We can now prove the main result Proposition~\ref{prop:analytic-finite-jet-recovery}:

\begin{proof}[Proof of Proposition~\ref{prop:analytic-finite-jet-recovery}]
We fix $0\leq j<d$ and use the following commutative diagram
with exact rows:
\[
\begin{tikzcd}[column sep=1.5em,row sep=large]
 0\arrow[r]
 &I_j\arrow[r,hook]\arrow[d,hook,"\epsilon^0"']
 &S_j^{q}\arrow[r,two heads]\arrow[d,hook,"\epsilon^0"']
 &\mathsf S(M_n)_j\arrow[r]\arrow[d,hook,"\epsilon^0"]&0\\
 0\arrow[r]
 &\bigl(I[\epsilon^{\pm1}]\bigr)_j\arrow[r,hook]
    \arrow[d,two heads,"\gr^j(\tau_d)"']
 &\bigl(S[\epsilon^{\pm1}]^{q}\bigr)_j\arrow[r,two heads]
    \arrow[d,two heads,"\gr^j(\tau_d)"']
 &\bigl(\mathsf S(M_n)[\epsilon^{\pm1}]\bigr)_j\arrow[r]
    \arrow[d,two heads,"\gr^j(\tau_d)"]&0\\
 0\arrow[r]
 &\gr_\omega^j\tau_d(K_n^\circ)\arrow[r,hook]
 &S_j^{q}\arrow[r,two heads]
 &\gr_{Q,n}^j\overline M_{n,d}^\circ\arrow[r]&0.
\end{tikzcd}
\]
The first two rows come from \eqref{eq:jet-bounded-presentation}
and its scalar extension to $S[\epsilon^{\pm1}]$; the last row comes from the
induced filtration on $\tau_d(K_n^\circ)$ and the quotient filtration.
All maps from the first row to the second row send $v$ to $\epsilon^0v$.

The composition in middle column $S_j^{q}\to(S[\epsilon^{\pm1}]^{q})_j\to S_j^{q}$ is the identity,
which is in particular injective. Thus the composition $I_j\to(I[\epsilon^{\pm1}])_j\to\gr_{\omega}^j\tau_d(K_n^\circ)$ is also injective. 

On the other hand, the map $(S[\epsilon^{\pm1}]^{q})_j\to S_j^{q}$
is given by
\[
 \bigl(S[\epsilon^{\pm1}]^{q}\bigr)_j
 =\bigoplus_{k\in\mathbb Z}S_{j-nk}^{q}\epsilon^k\xrightarrow{\gr^j(\tau_d)}S_j^{q},
 \qquad
 \gr^j(\tau_d)\Bigl(\sum_kv_k\epsilon^k\Bigr)=v_0.
\]
The subspace also has the decomposition \[(I[\epsilon^{\pm1}])_j=\bigoplus_{k\in\mathbb Z}I_{j-nk}\epsilon^k\]
Hence the surjective map $(I[\epsilon^{\pm1}])_j\twoheadrightarrow\gr_{\omega}^j\tau_d(K_n^\circ)$ factors through 
\(I_j\twoheadrightarrow\gr_{\omega}^j\tau_d(K_n^\circ).\)
We therefore obtain
\[
 I_j\xrightarrow{\sim}\gr_\omega^j\tau_d(K_n^\circ).
\]
Taking cokernels in the first and last rows then gives
\[
 \mathsf S(M_n)_j\xrightarrow{\sim}
 \gr_{Q,n}^j\overline M_{n,d}^\circ.
\]
The quotient filtration has zero graded pieces outside this range $0\le j<d$.
These maps are induced by the same graded presentation, so they give
\[
 \gr_{Q,n}^*\overline M_{n,d}^\circ
 \cong\frac{S^{q}}{I+(S^{q})_{\geq d}}
 \cong\frac{\mathsf S(M_n)}{\mathsf S(M_n)_{\geq d}},
\]
as desired.
\end{proof}

\section{\texorpdfstring{$D^\times$}{D-times}-equivariance and the \texorpdfstring{$J$}{J}-supported characteristic cycle}
\label{sec:equivariance}

\subsection{Filtered distribution algebras for \texorpdfstring{$U_D^1$}{UD1} and \texorpdfstring{$D^\times$}{D-times}}
\label{sec:D-equivariance}

In this subsection, we discuss the detailed properties of distribution algebras of the group $D^\times$. 

We first briefly recall the definitions and notions for $D$.
As a vector space, \[D=\mathbb{Q}_{p^2}\oplus\mathbb{Q}_{p^2}\varpi_D,\] with multiplication relations: \[ \varpi_D^2=p,\quad a\varpi_D=\varpi_D\sigma(a),\]
where $\sigma\in\mathrm{Gal}(\mathbb{Q}_{p^2}/\mathbb{Q}_p)$ is the unique non-trivial Frobenius element.
This $D$ admits a maximal order $\mathcal{O}_D=\mathbb{Z}_{p^2}\oplus\mathbb{Z}_{p^2}\varpi_D$. 
Then we have $D^\times=\cO_D^\times\rtimes\varpi_D^\mathbb{Z}$ with center $Z_D=\mathbb{Q}_p^\times$. There is a unique maximal pro-$p$ subgroup of $\cO_D^\times$, called $U_D^1=1+\varpi_D\cO_D$:
\[
 \cO_D^\times/U_D^1\cong H,\qquad \cO_D^\times=U_D^1\rtimes H,
\]
where $H$ is the Teichm\"uller subgroup $[\mathbb{F}_{p^2}^\times]$ fixed in the introduction.

Fix $\xi\in\F_{p^2}\backslash\F_p$. We use the standard ordered coordinates as in \cite{HuWangModP}:
\[
 \gamma_1=1+\varpi_D,
 \qquad
 \gamma_2=1+\varpi_D[\xi],
 \qquad
 \gamma_3=\gamma_1\gamma_2\gamma_1^{-1}\gamma_2^{-1},
 \qquad
 \gamma_4=1+p.
\]
As noted in \cite{HuWangModP} and \cite{Lazard1965}, $\gamma_1$ and $\gamma_2$ are minimal topological generators of $U_D^1/Z_D^1$, while $\gamma_1, \gamma_2, \gamma_3$ give the coordinates on $U_D^1/Z_D^1$ as a $p$-adic manifold.
The valuation on $U_D^1$ is defined by 
\[\omega(\gamma_1)=\omega(\gamma_2)=\frac{1}{2},\qquad\omega(\gamma_3)=\omega(\gamma_4)=1.\]

For a fixed unitary central character $\zeta:Z_D\to\cO_E^\times$, we define 
\begin{align*}
 \Lambda_\zeta^\circ
 &:=\cO_E\llbracket U_D^1\rrbracket/
   \overline{(z-\zeta(z):z\in Z_D^1)},\\
    D_{n,\zeta}
 &:=D_n(U_D^1,E)\widehat\otimes_{D_n(Z_D^1,E),\zeta}E,\\
 D_\zeta
 &:=D(U_D^1,E)\widehat\otimes_{D(Z_D^1,E),\zeta}E.
\end{align*}
We omit the group from the notation since we mainly work on the group $U_D^1$.

Since $Z_D^1$ is pro-$p$, the reduction of $\zeta|_{Z_D^1}$ modulo $\varpi_E$ is
trivial.  Consequently
\[
 \overline{\Lambda}_{\bar\zeta}:=\Lambda_\zeta^\circ/\varpi_E
 \cong k\llbracket U_D^1/Z_D^1\rrbracket.
\]
Set once and for all
\[
 S:=\gr_{\omega}^*\overline{\Lambda}_{\bar\zeta}
 \cong\gr_{\omega}^*k\llbracket U_D^1/Z_D^1\rrbracket.
\]

\begin{lemma}\label{omega=m-adic}
For $G=U_D^1/Z_D^1$ and integers $j\geq0$,
\[
 \mathrm{Fil}_\omega^{j/2}k\llbracket G\rrbracket=\mD^j,
 \qquad
 \gr_\omega^{j/2}k\llbracket G\rrbracket
       \xrightarrow{\sim}\gr_{\mD}^jk\llbracket G\rrbracket.
\]
\end{lemma}
\begin{proof}
This is \cite[Lemma~2.4]{HuWangModP}, after extending scalars to $k$.
\end{proof}

\paragraph{Full $D^\times$ representation}
We now define distribution algebras over the full group $D^\times$ with central character $\zeta$.
In each of the following rings, let $I_\zeta$ be the closed ideal generated by $\{z-\zeta(z)\}$ for all $z\in Z_D$. We define
\begin{align*}
 \Lambda_{D,\zeta}^\circ
 &\coloneqq(\cO_E\llbracket U_D^1\rrbracket\otimes_{\mathbb{Z}_p[U_D^1]}\mathbb{Z}_p[D^\times])/I_\zeta,\\
 \mathscr D_{n,\zeta}
 &\coloneqq(D_n(U_D^1,E)\otimes_{\mathbb{Z}_p[U_D^1]}\mathbb{Z}_p[D^\times])/I_\zeta,\\
 \overline{\Lambda}_{D,\bar\zeta}
 &\coloneqq\Lambda_{D,\zeta}^\circ/\varpi_E, \quad \mathscr{D}_\zeta\coloneqq(D(U_D^1,E)\otimes_{\mathbb{Z}_p[U_D^1]}\mathbb{Z}_p[D^\times])/I_\zeta.
\end{align*}

Let $\widetilde H=\langle H,\varpi_D\rangle$ and put $\Delta:=D^\times/(Z(D^\times)U_D^1)\cong (\mathbb{F}_{p^2}^\times/\mathbb{F}_p^\times)\rtimes\mathbb{Z}/2$. It has representatives
in $H\cup\varpi_DH$.  After central reduction, $\varpi_D^2=p$ acts by
the unit $\zeta(p)$.

Then we have a simple description of the rings above as a crossed product of rings for $U_D^1$ by the finite group $\Delta$:
The argument is exactly the same as that in \cite{Passman1989} and we omit it here.
\begin{proposition}\label{prop:equivariant-crossed-product-filtration}
      We have natural isomorphisms:
      \begin{align*}
 \Lambda_{D,\zeta}^\circ
 &=\Lambda_\zeta^\circ\rtimes_\zeta\Delta
   =\bigoplus_{\delta\in\Delta}\Lambda_\zeta^\circ[\delta],\\
 \mathscr D_{n,\zeta}
 &=D_{n,\zeta}\rtimes_\zeta\Delta
   =\bigoplus_{\delta\in\Delta}D_{n,\zeta}[\delta],\\
 \overline{\Lambda}_{D,\bar\zeta}
 &=\Lambda_{D,\zeta}^\circ/\varpi_E
   \cong\overline{\Lambda}_{\bar\zeta}\rtimes_{\zeta}\Delta=\bigoplus_{\delta\in\Delta}\overline{\Lambda}_{\bar\zeta}[\delta].
\end{align*}
Here $\rtimes_\zeta$ denotes the finite twisted crossed product for
conjugation by $\widetilde H$, with the relations
$[z]=\zeta(z)$ for $z\in Z_D$.

Therefore, we have an isomorphism of Fr\'echet-Stein algebras
\[\mathscr{D}_\zeta\cong\varprojlim_n\mathscr{D}_{n,\zeta}.\]
\end{proposition}

For the distributions over $D^\times$, we have 
\[
 \Lambda_\zeta=\Lambda_\zeta^\circ[1/\varpi_E],\qquad
 \Lambda_{D,\zeta}=\Lambda_{D,\zeta}^\circ[1/\varpi_E],\qquad
 \mathscr D_\zeta=D_\zeta\rtimes_\zeta\Delta.
\]

The conjugation action of $D^\times$ preserves the valuation $\omega$ on $U_D^1$,
the coordinate description above makes the central reductions strict,
and
\(D^\times/(Z_DU_D^1)=\Delta\)
has order prime to $p$.  Therefore we define the corresponding filtrations on these three rings $R$ above by the component-wise filtration with each $[\delta]$ has filtration degree $0$, i.e.
\[\mathrm{Fil}_{n}^*(R\rtimes_\zeta\Delta)\coloneqq\bigoplus_{\delta\in\Delta}\mathrm{Fil}_n^*R\cdot[\delta].\]

These filtrations are complete. The following lemma follows immediately from the definition:

\begin{lemma}\label{lem:equivariant-module-filtration}
For all sufficiently large \(n\),
\[
 \gr_n^*\Lambda_{D,\zeta}^\circ
 \cong S[\epsilon]\rtimes_\zeta\Delta,
 \qquad
 \gr_n^*\mathscr D_{n,\zeta}
 \cong S[\epsilon^{\pm1}]\rtimes_\zeta\Delta,
\]
and hence
\[
 \gr_\epsilon^0\gr_n^*\mathscr D_{n,\zeta}
 \cong S\rtimes_\zeta\Delta.
\]
All filtrations on filtered \(\mathscr D_{n,\zeta}\) modules over are \(\Delta\)-stable.
\end{lemma}

We now choose $H$-eigenvectors in $\overline\Lambda_{\bar\zeta}$.
Following \cite[\S2.3]{HuWangModP}, put
\[Y=\sum_{\lambda\in\mathbb{F}_{p^2}^\times}\lambda^{-1}(1+\varpi_D[\lambda]), Z=\sum_{\lambda\in\mathbb{F}_{p^2}^\times}\lambda^{-p}(1+\varpi_D[\lambda])\]
By \cite[Lemma~2.7]{HuWangModP}, $Y,Z$ generate the maximal
ideal $\mD$ of $\overline\Lambda_{\bar\zeta}$.
Write $y,z$ for their principal symbols in $S$.
Their $H$-characters are $\alpha$ and $\alpha^{-1}$, where
$\alpha(h)=h^{p-1}$, and
\[
 \varpi_DY\varpi_D^{-1}=Z,
 \qquad
 \varpi_DZ\varpi_D^{-1}=Y.
\]

In these coordinates, \cite[\S2.3]{HuWangPadic} identifies
\[
 S\cong k\langle y,z,h\rangle/
       (yz-zy-h,\ yh-hy,\ zh-hz),\qquad
 \deg_\omega y=\deg_\omega z=\tfrac12,\quad\deg_\omega h=1.
\]
In particular $h$ is central and regular, $S/(h)=k[y,z]$, and the
ordered monomials in the three generators form a PBW basis.

\subsection{The \texorpdfstring{$J$}{J}-supported characteristic cycle}
\label{sec:quaternionic}
We apply Subsection~\ref{sec:uniform} to the ring $S$ above.
For a fixed sufficiently large $n$, we have
\[
 \gr_n^*D_{n,\zeta}=S[\epsilon^{\pm1}],\qquad
 \gr_n^*\Lambda_\zeta^\circ=S[\epsilon].
\]
The PBW filtration shows that $S$ is noetherian. The filtration
on $D_{n,\zeta}$ is complete and quasi-integral, with noetherian
associated graded ring, so it is Zariskian. We use good module
filtrations with rational generator degrees.

Following \cite[\S3]{HuWangPadic}, we put
\[
 J=(yz,zy)=(h,yz)\subset S,\qquad
 A_0=S/J=k[y,z]/(yz).
\]
The commutator ideal of $S$ is $(h)$, which is contained in $J$.
Thus $S/J$ is commutative and satisfies the condition in Subsection~\ref{sec:uniform}.
Its two branches are $V(z)$ and $V(y)$.

Let $N$ be a finite graded $S$-module with $J^aN=0$.
We use Hu--Wang's branch multiplicities
\cite[\S3.2, equation (3.2)]{HuWangPadic}. Since
\[
 (A_0)_{(z)}=k(y),\qquad (A_0)_{(y)}=k(z),
\]
the one-dimensional cycle \eqref{eq:q-cycle} is defined by
\begin{align}
 \CC_J(N)
 &=\sum_{i=0}^{a-1}
   \dim_{k(y)}\bigl((J^iN/J^{i+1}N)\otimes_{A_0}k(y)\bigr)[V(z)]
 \label{eq:D-cycle}\\
 &\quad+\sum_{i=0}^{a-1}
   \dim_{k(z)}\bigl((J^iN/J^{i+1}N)\otimes_{A_0}k(z)\bigr)[V(y)].
 \nonumber
\end{align}
As in \cite{HuWangPadic}, $\mu(N)$ is the sum of the two coefficients.
By Lemma~\ref{lem:local-length}, both $\CC_J$ and $\mu$ are
additive in exact sequences and unchanged by grading shifts.
We use the one-dimensional cycle throughout this subsection;
in particular, it is zero for finite-dimensional graded modules.

We now apply the general results to show that this cycle is well-defined.

\begin{corollary}\label{thm:quaternionic-CC}
Let $M$ be a finite $D_{n,\zeta}$-module. Suppose one residual
symbol $\mathsf S(M)$, obtained from a good filtration and an
integral graded lattice, is annihilated by a power of $J$.
Then every such residual symbol is annihilated by a power of
$J$, and
\[
 \CC_J(\mathsf S(M)),\qquad \mu(\mathsf S(M))
\]
are independent of the good filtration, integral graded lattice,
and filtered presentation. These modules form a Serre subcategory,
and the cycle and multiplicity are additive in exact sequences.
\end{corollary}
\begin{proof}
We apply Theorems~\ref{thm:main} and \ref{thm:additivity} with
$A=D_{n,\zeta}$, $J=(h,yz)$, and $q=1$. The ring $A_0$ has
dimension one, and \eqref{eq:D-cycle} gives
\[
 \CC_J(\mathsf S(M))=\CC_1(M).
\]
This proves all the assertions.
\end{proof}

Thus, by abuse of notations, we omit the symbol $\mathsf{S}$ and write, for finite $D_{n,\zeta}$ modules $M$, 
\[\CC_J(M)\coloneqq\CC_J(\mathsf S(M)),\qquad\mu(M)\coloneqq\mu(\mathsf{S}(M))\]
for some choice of good filtration and integral graded lattice. 

We first identify the modules whose cycle is zero.

\begin{corollary}\label{cor:D-zero}
For a finite graded $J$-power-torsion $S$-module $N$,
$\CC_J(N)=0$ if and only if $\dim_kN<\infty$.
Consequently, for $M$ in Corollary~\ref{thm:quaternionic-CC},
$\CC_J(M)=0$ implies $\dim_EM<\infty$.
\end{corollary}
\begin{proof}
Since $A_0=k[y,z]/(yz)$ has dimension one and $\gr_JN$ is graded,
\[
 \CC_J(N)=0
 \quad\Longleftrightarrow\quad
 \Supp_{A_0}(\gr_JN)\subseteq\{(y,z)\}
 \quad\Longleftrightarrow\quad \dim_kN<\infty.
\]
If $N=\mathsf S(M)=L/\epsilon L$ is finite-dimensional, we lift a homogeneous
basis and use $\epsilon$-regularity to see that $L$ is finite graded-free
over $k[\epsilon]$. Therefore by \cite[Theorem~I.4.2.4.(5)]{LvO96} we get
\[
 \dim_EM=\rank_{k[\epsilon^{\pm1}]}\gr_FM=\dim_kN<\infty.
\]
\end{proof}

We now use the group action to compare the two branch coefficients.

\begin{corollary}
\label{cor:D-symmetry}
Let $M$ be a finite $\mathscr D_{n,\zeta}$-module whose residual
symbol is annihilated by a power of $J$. The coefficients of
$[V(z)]$ and $[V(y)]$ in $\CC_J(M)$ are equal.
Consequently, $\mu(M)$ is an even
nonnegative integer and
\[
 \CC_J\bigl(M\bigr)=\frac{\mu(M)}{2}
                    \bigl([V(z)]+[V(y)]\bigr).
\]
\end{corollary}
\begin{proof}
By Lemma~\ref{lem:equivariant-module-filtration}, we choose a
$\widetilde H$-stable filtration and integral graded lattice. Conjugation
by $\varpi_D$ interchanges $y,z$ and sends $h$ to $-h$.

It preserves $J$ and hence acts on each quotient of the
$J$-filtration of $\mathsf S(M)$, interchanging its localizations
at $(z)$ and $(y)$ in $A_0$. These identifications preserve the
dimensions over the corresponding fields $k(y)$ and $k(z)$.

We sum over the $J$-filtration in \eqref{eq:D-cycle} and apply
Corollary~\ref{thm:quaternionic-CC} to show that the result is independent
of the stable filtration and lattice. This also agrees with \cite[Lemma~3.6]{HuWangPadic}.
\end{proof}

\subsection{Purity across filtrations}
\label{sec:nonexistence}
Throughout this subsection, we fix a radius index $n$ and consider
finite $D_{n,\zeta}$-modules $M$ whose residual symbols are
annihilated by a power of $J$.
We aim to construct a quotient $Q$ of $M$ with a good filtration and an
integral graded lattice such that
\[
 \mathsf S(Q)\text{ is pure of grade }2,\qquad
 \CC_J(Q)=\CC_J(M),
\]
provided that the common cycle is nonzero.

\paragraph{Grades and dimensions}
\label{sec:grade-comparison}

We write $\operatorname{GKdim}_R M$ for Gelfand--Kirillov dimension, or equivalently the canonical dimension 
in the dimension theory of Auslander regular rings. 
Recall that \[S=\gr_\omega k[\![U_D^1/Z_D^1]\!]\] is an Auslander regular ring. And so are $S[\epsilon]$ and $S[\epsilon^{\pm1}]$. 

By filtered grade comparison
\cite[Chapter~III, Theorem~2.5.2]{LvO96}, for every nonzero finite $R$-module $M$,
we get 
\[
 j_R(M)+\operatorname{GKdim}_R M
 =\operatorname{injdim}_R R
 =\begin{cases}
 3,&R=S,\\
 4,&R=S[\epsilon]\text{ or }S[\epsilon^{\pm1}].
 \end{cases}
\]
Thus Gelfand--Kirillov dimension agrees with canonical dimension
for these rings.

Let $V$ be a nonzero finite $S[\epsilon]$-module on which
$\epsilon$ is regular. Since $\epsilon$ is central, localization
preserves Gelfand--Kirillov dimension. Hence
\[
 \operatorname{GKdim}_{S[\epsilon]}V
 =\operatorname{GKdim}_{S[\epsilon^{\pm1}]}
      V[\epsilon^{-1}],\qquad
 j_{S[\epsilon]}(V)
 =j_{S[\epsilon^{\pm1}]}(V[\epsilon^{-1}]).
\]
We now compute the grade over the rings
$S/(h)=k[y,z]$ and $S/J=k[y,z]/(yz)$.

\begin{lemma}\label{lem:residual-grade}
Let $I\in\{(h),J\}$ and let $N\ne0$ be a finite $S$-module
with $I^aN=0$ for some $a$. Then
\[
 j_S(N)+\dim\Supp_{S/I}(\gr_I N)=3.
\]
For $I=J$, this support has dimension at most one.
\end{lemma}
\begin{proof}
Each nonzero $U_i=I^iN/I^{i+1}N$ is killed by $h$.
By the Rees lemma \cite[Proposition~3.4(b)]{Levasseur1992}, we get
\[
 j_S(U_i)=1+j_{k[y,z]}(U_i)
 =3-\dim\Supp_{S/I}U_i.
\]
Using this composition series, we get
\[
 j_S(N)=\min_{U_i\ne0}j_S(U_i)
 =3-\dim\Supp_{S/I}(\gr_I N).
\]
Finally, when $I=J$, the ring $S/J=k[y,z]/(yz)$ has dimension one.
\end{proof}

We now compare the grades of a distribution module and its residual symbol.

\begin{corollary}\label{cor:quaternionic-grade}
For every nonzero finite $D_{n,\zeta}$-module $V$, we have
\[
 j_{D_{n,\zeta}}(V)=j_S(\mathsf S(V)).
\]
If its residual symbol is annihilated by a power of $J$, then
\[
 j_{D_{n,\zeta}}(V)=
 \begin{cases}
 2,&\CC_J(\mathsf S(V))\ne0,\\
 3,&\CC_J(\mathsf S(V))=0.
 \end{cases}
\]
\end{corollary}\begin{proof}
Choose a good filtration $F$ and an integral graded lattice
$L\subset\gr_FV$. Since $\epsilon$ is homogeneous of positive
degree and $L$ is $\epsilon$-torsion-free,
\cite[Corollary~4.4]{Levasseur1992} tells us
\[
 j_{S[\epsilon]}(L/\epsilon L)=j_{S[\epsilon]}(L)+1.
\]
On the other hand, by the Rees lemma \cite[Proposition~3.4(b)]{Levasseur1992},
we have
\[
 j_S(L/\epsilon L)
 =j_{S[\epsilon]}(L/\epsilon L)-1
 =j_{S[\epsilon]}(L).
\]
Thus, by the localization comparison above and filtered grade comparison
\cite[Chapter~III, Theorem~2.5.2]{LvO96}, we can identify
\[
 \begin{aligned}
 j_{D_{n,\zeta}}(V)
 &=j_{S[\epsilon^{\pm1}]}(\gr_FV)\\
 &=j_{S[\epsilon]}(L)
 =j_S(L/\epsilon L)
 =j_S(\mathsf S(V)).
 \end{aligned}
\]
The assertion about the cycle follows from Lemma~\ref{lem:residual-grade} and Corollary~\ref{cor:D-zero}.
\end{proof}

\paragraph{Pure residual symbols}

We use the notation $\Delta_R^i$ and the descending grade indexing
of \cite[\S8]{ST03} for the canonical dimension filtration.
We recall its construction from the biduality spectral sequence;
see also \cite[\S3.1]{Venjakob2002}.
Let $R$ be Auslander regular and $V$ a finite left $R$-module, and put
\[
 d_R=\operatorname{injdim}_R R.
\]
Since $V$ has a finite resolution by finite projective modules, we have
\[
 V\xrightarrow{\sim}
 \mathbf R\!\operatorname{Hom}_{R^{\mathrm{op}}}
 \bigl(\mathbf R\!\operatorname{Hom}_R(V,R),R\bigr).
\]
The associated spectral sequence, with $a\geq0$ and $b\leq0$, is given by
\[
 E_2^{a,b}=
 \operatorname{Ext}_{R^{\mathrm{op}}}^{a}
 \bigl(\operatorname{Ext}_R^{-b}(V,R),R\bigr)
 \quad\Longrightarrow\quad
 V[0].
\]
It gives us the canonical dimension filtration
\[
 \begin{gathered}
 V=\Delta_R^0(V)\supseteq\Delta_R^1(V)\supseteq\cdots
   \supseteq\Delta_R^{d_R+1}(V)=0,\\
 \Delta_R^i(V)/\Delta_R^{i+1}(V)\cong E_\infty^{i,-i}.
 \end{gathered}
\]
By the Auslander condition, each nonzero quotient is pure of grade $i$.
We also have the following descriptions by
\cite[Propositions~3.5--3.6]{Venjakob2002}:
\[
 \Delta_R^i(V)=\sum_{\substack{U\subseteq V,\ j_R(U)\geq i}}U,
 \qquad \Delta_R^i\bigl(V/\Delta_R^i(V)\bigr)=0.
\]
In particular,
\[
 j_R(V)=a
 \quad\Rightarrow\quad
 V/\Delta_R^{a+1}(V)\text{ is pure of grade }a.
\]
By functoriality of the spectral sequence, compatible semilinear actions
preserve the filtration.

We now apply Gabber's maximality principle to the $\epsilon$-filtration; see also \cite[\S8.5]{ArdakovWadsley2013}.
We first recall the statement below. For the proof, see \cite[Theorem~1.14]{Bjork1989}.
\begin{theorem}
      \label{thm:Gabber-max}
      Let $R$ be an Auslander-Gorenstein ring and let $X\subset Y$ be two $R$-modules such that $X$ is a finitely generated pure $R$-module, and each finitely generated submodule of $Y$ is pure. 
      Then $Y$ contains a unique largest finitely generated submodule $Z$ with $X\subset Z$ such that 
      \[j_{R}(Z/X)\ge j_R(X)+2. \] 
\end{theorem}

\begin{lemma}\label{lem:pure-lattice}
Let $P$ be a finite graded $S[\epsilon^{\pm1}]$-module pure of
grade $2$, and suppose one lattice reduction is annihilated by
a power of $J$. For every integral graded lattice $L\subset P$,
there is a largest finite $S[\epsilon]$-submodule $L^\dagger\subset P$
such that
\[
 L\subseteq L^\dagger,\qquad \dim_k(L^\dagger/L)<\infty.
\]
It is an integral graded lattice, and $L^\dagger/\epsilon L^\dagger$
is pure of grade $2$, is annihilated by a power of $J$, and has the
same cycle as $L/\epsilon L$.
If $L$ is $\widetilde H$-stable, then so is $L^\dagger$.
For a compatible $\widetilde H$-action on $P$, such a stable $L$ exists.
\end{lemma}
\begin{proof}
For every nonzero finite $S[\epsilon]$-submodule $V\subset P$,
$\epsilon$ is regular on $V$. The localization comparison above tells us
\[
 j_{S[\epsilon]}(V)
 =j_{S[\epsilon^{\pm1}]}(V[\epsilon^{-1}])=2.
\]
Thus every such submodule is pure of grade $2$.
Then the Gabber's maximality principle (Theorem~\ref{thm:Gabber-max})
  gives us a largest finite extension $L^\dagger\subset P$ of $L$ with
\[
 j_{S[\epsilon]}(L^\dagger/L)\geq4=\dim S[\epsilon]
 \quad\Longleftrightarrow\quad \dim_k(L^\dagger/L)<\infty.
\]
The span of the homogeneous components of $L^\dagger/L$ is again
a finite-dimensional submodule of $P/L$. By maximality,
$L^\dagger$ is graded. What's more, $L^\dagger$ is an integral graded lattice since $L^\dagger\supset L$.

Put $U=\Delta_S^3(L^\dagger/\epsilon L^\dagger)$ and let
$W\subseteq L^\dagger$ be its inverse image. Then by definition $\epsilon L^\dagger\subset W$, and thus
\[
 L^\dagger\subseteq\epsilon^{-1}W\subset P,\qquad
 (\epsilon^{-1}W/L^\dagger)\langle-n\rangle\cong U.
\]
The grade--dimension formula gives $j_S(U)\ge3\Rightarrow\dim_k U<\infty$. Hence
\[
 \dim_k(\epsilon^{-1}W/L)<\infty
 \quad\Rightarrow\quad\epsilon^{-1}W=L^\dagger
 \quad\Rightarrow\quad U=0.
\]
The reduction $L^{\dagger}/\epsilon L^{\dagger}$ is nonzero by the graded Nakayama
\cite[Lemma~6.2]{NegronPevtsova2023}. Applying
Proposition~\ref{prop:lattice-invariance}, we get the same characteristic cycle
and annihilation by a power of $J$. By the Lemma~\ref{lem:residual-grade}
and the fact that $U=0$, we can deduce that
\[
 0\ne V\subseteq L^\dagger/\epsilon L^\dagger
 \quad\Rightarrow\quad 2\leq j_S(V)<3
 \quad\Rightarrow\quad j_S(V)=2.
\]As a consequence, $L^{\dagger}/\epsilon L^{\dagger}$ is pure. 

Finally, if $L$ is stable under $\widetilde{H}$, uniqueness gives $gL^\dagger=L^\dagger$.
We obtain a stable $L$ by summing its translates under the finite
crossed product in Proposition~\ref{prop:equivariant-crossed-product-filtration}.
\end{proof}

We can now construct the required quotient of the distribution module.

\begin{proposition}
\label{prop:analytic-purification}
Let $M$ be a finite $\mathscr D_{n,\zeta}$-module.
There is a finite-dimensional $\mathscr D_{n,\zeta}$-submodule
$T\subset M$ such that $Q=M/T$ satisfies
\[
 \CC_J(Q)=\CC_J(M),\qquad
 \mu(Q)=\mu(M).
\]
If $\mu(M)=0$, then $Q=0$. Otherwise, $Q$ admits
a compatible good filtration and a stable graded lattice such that the corresponding
residual symbol $\mathsf S(Q)$ is pure of grade $2$.
\end{proposition}
\begin{proof}
We take the following term from the canonical dimension filtration and note that $\mathscr{D}_{n, \zeta}$ is finite over $D_{n,\zeta}$,
\[
 T=\Delta_{D_{n,\zeta}}^3(M)=\Delta_{\mathscr{D}_{n,\zeta}}^3(M),
 \qquad Q=M/T.
\]
Since $T$ is a $\mathscr{D}_{n, \zeta}$-module, it carries the full $D^\times$-action, and so does $Q$. By the Corollaries~\ref{cor:quaternionic-grade} and \ref{cor:D-zero} above, we get
\[
 j_{D_{n,\zeta}}(T)\geq3
 \quad\Rightarrow\quad \CC_J(T)=0
 \quad\Rightarrow\quad \dim_E T<\infty.
\]
Applying Corollary~\ref{thm:quaternionic-CC} to the exact sequence
\[
 0\longrightarrow T\longrightarrow M\longrightarrow Q=M/\Delta^3_{D_{n,\zeta}}(M)
 \longrightarrow0
\]
we get the asserted cycle and multiplicity equalities.
If $\mu(M)=0$, Corollary~\ref{cor:quaternionic-grade}
gives us $j_{D_{n,\zeta}}(M)\geq3$, hence $T=M$.
Otherwise, the same corollary and the purity of the canonical
filtration quotients give us
\[
 j_{D_{n,\zeta}}(M)=2
 \quad\Rightarrow\quad Q\text{ is pure of grade }2.
\]

Choose a $\widetilde H$-stable good filtration on $Q$ with
rational generator degrees.
We now apply Theorem~\ref{thm:gabber-filtration} to $Q$ and this filtration and with
\[
 A=\mathscr D_{n,\zeta},\qquad \gr A=S[\epsilon^{\pm1}]\rtimes_{\bar{\zeta}}\Delta,\qquad j=2.
\]
This gives us a Gabber filtration $F$ on $Q$ that is stable under the $\widetilde{H}$-action and for which $\gr_FQ$ is pure of grade $2$.
Furthermore, Lemma~\ref{lem:pure-lattice} gives us a stable lattice $L$
whose reduction $L/\epsilon L$ is pure of grade $2$, with the required cycle and multiplicity.
\end{proof}

\subsection{Representations with multiplicity 2}

In this subsection, we prove the non-existence of representations of multiplicity two. 
This is similar to \cite[Corollary~3.12]{HuWangPadic} but in the locally analytic setting.

We first classify pure residual symbol modules with multiplicity $2$.
\begin{corollary}\label{cor:two-symbol-classification}
For a finite pure graded $S\rtimes_{\bar\zeta}\Delta$-module $N$ 
annihilated by a power of $J$ and satisfying $\mu(N)=2$, there are an
$H$-character $\chi$ and a common grading shift $\ell$ such that
\[
 N\cong\bigl((\chi\otimes k[y])\oplus
       (\chi^\sigma\otimes k[z])\bigr)\langle-\ell\rangle
\]
or
\[
 N\cong(\chi\otimes k[y,z]/(yz))\langle-\ell\rangle,
 \qquad\chi=\chi^\sigma.
\]
In particular, $JN=0$ and $N$ is generated in one degree over $S$.
\end{corollary}
\begin{proof}
Being an $S\rtimes_{\bar\zeta}\Delta$ module is equivalent to be a $S=\gr_\omega\overline{\Lambda}_{\bar\zeta}$ module with a compatibile $\widetilde{H}$-action. 
By assumption $N$ is pure and annihilated by a power of $J$. Thus the assertion follows from \cite[Lemma~3.7]{HuWangPadic}.
\end{proof}

Now let $L_n$ be a finite $\mathscr D_{n,\zeta}$-module with a pure
residual symbol of multiplicity two, as obtained in
Proposition~\ref{prop:analytic-purification}.
After a common grading shift, we assume that the module $\mathsf S(L_n)$ is generated in degree zero.
As in the final steps of the proof of Lemma~\ref{lem:two-stage-strictness}, 
we choose the generators of $L_n$ such that they and their principal symbols 
define the same radius filtration and the chosen
$\epsilon$-filtration. Thus Lemma~\ref{lem:analytic-jet-relations}
also applies to the lattice in Lemma~\ref{lem:pure-lattice}.

We use the notation and construction of
Subsection~\ref{sec:analytic-to-residual} with $G=U_D^1$ and
$M_n=L_n$. For $d>0$ and $n\geq1+d+1/e$, we get
\[
 \tau_d:D_n^\circ\twoheadrightarrow
       k\llbracket U_D^1\rrbracket/J_d,\qquad
 \overline L_{n,d}^\circ=L_n^\circ/\mathfrak a_{n,d}L_n^\circ,
 \qquad L_n^\circ=\mathrm{Fil}_n^0L_n,
\]
with the quotient filtration used there.
Conjugation preserves the filtrations and
$\tau_d(gag^{-1})=g\tau_d(a)g^{-1}$, so the construction is
$D^\times$-equivariant. Every resulting quotient satisfies
\[
 \varpi_E\overline{L}_{n,d}^{\circ}=0,\qquad
 (z-1)\overline{L}_{n,d}^{\circ}=0\quad(z\in Z_D^1),
\]
and is a finite smooth $D^\times$-representation with central
character $\bar\zeta$.

To apply the construction over $U_D^1$, we put
$\widetilde S:=\gr_\omega^*k\llbracket U_D^1\rrbracket$.
Since $\zeta(\gamma_4)-1\in\varpi_E\cO_E$, for $n>1$ we have
\[
 \sigma_n\bigl(\gamma_4-\zeta(\gamma_4)\bigr)
 =\sigma_n(\gamma_4-1)=:w_4,\qquad
 \deg w_4=1,\qquad \widetilde S/(w_4)\cong S.
\]
Thus the residual symbols over $\widetilde S$ have the additional
relation $w_4=0$.

\begin{lemma}\label{lem:two-symbol-truncation}
Let $L_n$ be as above, with $\mathsf S(L_n)$ generated in degree
zero. Let $d\in\tfrac12\mathbb Z_{>0}$ and suppose
$n\geq1+d+1/e$. Then
\[
 \gr_{Q,n}^{j/2}\overline L_{n,d}^{\circ}
 =\gr_{\mD}^{j}\overline L_{n,d}^{\circ}\qquad(j\geq0),
\]
compatibly with the $D^\times$-action. According to the two cases of
Corollary~\ref{cor:two-symbol-classification}, we have
\[
 \gr_{\mD}^*\overline L_{n,d}^\circ\cong
 \begin{cases}
 (\chi_n\otimes k[y]/(y^{2d}))\oplus
 (\chi_n^\sigma\otimes k[z]/(z^{2d})),&\text{case (1)},\\
 \chi_n\otimes k[y,z]/(yz,y^{2d},z^{2d}),&\text{case (2)}.
 \end{cases}
\]
Here $\chi_n$ is the character in that corollary, with
$\chi_n=\chi_n^\sigma$ in case~\textup{(2)}, and the degrees on the
right are $\deg y=\deg z=1$ in the $\mathfrak{m}$-adic filtration.
\end{lemma}
\begin{proof}
Before regrading, the generators have degree zero and the
relations over $\widetilde S$ are generated by
\[
 \begin{cases}
 z\mathbf e_1,\ h\mathbf e_1,\ w_4\mathbf e_1,
 \ y\mathbf e_2,\ h\mathbf e_2,\ w_4\mathbf e_2,&\text{case (1)},\\
 h,\ yz,\ w_4,&\text{case (2)}.
 \end{cases}
\]
Since
\[
 \deg y=\deg z=\tfrac12,\qquad \deg h=\deg w_4=1,
\]
all relations have degree at most one, independently of $n$. We take
$d'=1$ in Proposition~\ref{prop:analytic-finite-jet-recovery} and get
\[
 \gr_{Q,n}^*\overline L_{n,d}^\circ
 \cong\mathsf S(L_n)/\mathsf S(L_n)_{\geq d}.
\]
The central quotient is strict in the ordered coordinates. Hence the
degree-zero presentation and Lemma~\ref{omega=m-adic} give us
\[
 \mathrm{Fil}_{Q,n}^{j/2}\overline L_{n,d}^{\circ}
 =\bigl(\mathrm{Fil}_\omega^{j/2}\overline\Lambda_{\bar\zeta}\bigr)
   \overline L_{n,d}^{\circ}
 =\mD^j\overline L_{n,d}^{\circ}.
\]
Taking successive quotients gives us the first assertion. We regrade by
$j/2\mapsto j$. Since $2d$ is an integer, the truncation kills
$y^{2d}$ and $z^{2d}$. This gives us the two displayed modules.
\end{proof}

We now recall the finite-dimensional obstruction in \cite{HuWangPadic} used below.

\begin{theorem}[{\cite[Theorem~2.33]{HuWangPadic}}]
\label{thm:HW-finite-obstruction}
Let $p\geq5$ and let $\chi:\cO_D^\times\to k^\times$ be a
smooth character. There is no smooth $D^\times$-representation
$\pi$ over $k$ with a central character for which
\[
 \gr_{\mD}^*(\pi^\vee)\cong
 \bigl(\chi^{-1}\otimes k[y]/(y^{p+3})\bigr)
 \oplus
 \bigl((\chi^\sigma)^{-1}\otimes k[z]/(z^{p+3})\bigr)
\]
as graded $S$-modules with a compatible $H$-action.
Here $S$ is regraded by $s\mapsto2s$, so $\deg y=\deg z=1$,
and both generators have degree zero.
\end{theorem}

In this statement, $Z_D^1$ acts trivially, as it is pro-$p$ and
the central character takes values in $k^\times$.
Thus $\pi^\vee$ is a module over
$k\llbracket U_D^1/Z_D^1\rrbracket$, with its $\mD$-adic
filtration. 

We now apply this obstruction to the finite quotients constructed above.

\begin{theorem}\label{thm:no-two}
For all sufficiently large $n$ satisfying
\[
 n\geq1+\frac{p+4}{2}+\frac1e,
\]
there is no finite $\mathscr D_{n,\zeta}$-module $M$ such that
$\mathsf S(M)$ is annihilated by a power of $J$ and
$\mu(M)=2$.
\end{theorem}

\begin{proof}
Suppose such an $M_n$ exists. By
Proposition~\ref{prop:analytic-purification}, we take a quotient
$L_n$ with the same cycle and a pure residual symbol. We shift the
filtration so that $\mathsf S(L_n)$ is generated in degree zero.

By Corollary~\ref{cor:two-symbol-classification}, we make the following choices in the two cases:
\[\begin{cases}
 d=\frac{p+3}{2},\qquad
 V=\overline L_{n,d}^\circ,\qquad &\theta=\chi_n, \quad\text{case (1)};\\
 d=\frac{p+4}{2},\qquad
 V=\mD\overline L_{n,d}^\circ,\qquad &\theta=\chi_n\alpha,\quad \text{case (2)}.
\end{cases}
\]
Here $\alpha$ is the $H$-character of $y$. In the second case,
\[
 \mD^jV=\mD^{j+1}\overline{L}_{n,d}^{\circ},\qquad
 (\chi_n\alpha)^\sigma=\chi_n\alpha^{-1},
\]
since $\chi_n=\chi_n^\sigma$.

By assumption $n\geq1+d+1/e$ in both
cases. Thus, by Lemma~\ref{lem:two-symbol-truncation}, we get
\[
 \gr_{\mD}^*V
 \cong\bigl(\theta\otimes k[y]/(y^{p+3})\bigr)
       \oplus\bigl(\theta^\sigma\otimes k[z]/(z^{p+3})\bigr).
\]
The module $V$ is a finite smooth $D^\times$-representation
with central character $\bar\zeta$. We set
\[
 \pi=V^\vee,\qquad \chi=\theta^{-1},
\]
extending $\chi$ trivially across $U_D^1$.
Then $\pi$ is smooth with central character $\bar\zeta^{-1}$,
and $\pi^\vee=V$. The displayed isomorphism is compatible with
$H$, and both generators have degree zero. This contradicts
Theorem~\ref{thm:HW-finite-obstruction}.
\end{proof}

Since a locally analytic representation involves all radii, we get the following corollary.     
\begin{corollary}\label{cor:no-eventual-two}
There is no admissible locally analytic $D^\times$-representation
$\Pi^{\mathrm{an}}$ with central character $\zeta^{-1}$ such that,
for all sufficiently large $n$, the stages
\[
 M_n=\mathscr D_{n,\zeta}\otimes_{\mathscr D_\zeta}
           (\Pi^{\mathrm{an}})'_b
\]
have residual symbols annihilated by a power of $J$ and
$\mu(M_n)=2$.
\end{corollary}
\begin{proof}
We choose $n$ large enough that the two conditions hold and
the bound in Theorem~\ref{thm:no-two} is satisfied.
That theorem gives us a contradiction.
\end{proof}

\subsection{Residual symbols under infinitesimal characters with small multiplicity}
\label{sec:infinitesimal-geometric-two}
In general, the $J$-power annihilation of a residual symbol is hard to detect. 
In some cases, we can detect $J$-power annihilation using information from the representation, and these cases are enough for our application.

We now give a sufficient condition for $J$-power annihilation using an
infinitesimal character. Recall that $e=v_E(p)$. We work with the cofinal
radius indices
\[
 n=\frac{p^m}{e}
\]
with $m$ sufficiently large. Since $p\geq5$, we have
\[
 \left(1-\frac1p\right)p^m<ne=p^m<\frac{p-1}{2}p^m.
\]

We first compute the symbols of logarithms at these radii, following
\cite[Lemma~6.5]{Sorgdrager2026}.

\begin{lemma}\label{lem:cofinal-logarithm-symbol}
Let $g$ be an element of an ordered basis, with
$\omega(g)\in\{\tfrac12,1\}$. In $D_n(U_D^1,E)$, at the radii above, we have
\[
 \sigma_n(\log g)=\sigma_n(p)^{-m}\sigma_n(g-1)^{p^m},
 \qquad
 \deg_n(\log g)=p^m\omega(g)-ne m.
\]
\end{lemma}
\begin{proof}
Put $b=g-1$. In the convergent expansion
\[
 \log g=\sum_{r\geq1}\frac{(-1)^{r+1}}{r}b^r,
\]
the term with exponent $r$ has degree
$r\omega(g)-ne v_p(r)$. Among exponents with $v_p(r)=s$, the
least degree occurs at $r=p^s$. The least degree occurs uniquely at the exponent $p^m$
precisely when
\[
 \omega(g)(p-1)p^{m-1}<ne<\omega(g)(p-1)p^m.
\]
The chosen bounds on $ne$ imply both inequalities for
$\omega(g)=\tfrac12$ and $1$. The coefficient of $b^{p^m}$ is
$p^{-m}$, since $p$ is odd. This proves the assertion.
\end{proof}

We next apply this calculation to the Casimir element.
Choose a trace-zero unit $\vartheta\in\mathbb Z_{p^2}$, and write
$a=\vartheta^2\in\mathbb Z_p^\times$. In the Lie algebra of $D^\times$,
put
\[
 X_1=\varpi_D,\qquad X_2=\varpi_D\vartheta,\qquad X_3=p\vartheta.
\]
The elements $\exp(X_1),\exp(X_2),\exp(X_3),\exp(p)$ form an
ordered basis of $U_D^1$, with weights $\tfrac12,\tfrac12,1,1$.
Indeed, their logarithms form a $\mathbb Z_p$-basis of
$\varpi_D\cO_D$; see also the exponential coordinates in
\cite[\S6.5]{Sorgdrager2026}.

These coordinates give the same radius filtration by a simple computation; cf.\ \cite[Lemma~4.1]{ST03}.
Let $x_i=\sigma_n(\exp(X_i)-1)$. By direct computation in the variables $y, z, h$, we get  
\[x_1=-(y+z), \qquad x_2=\bar \vartheta(y-z), \qquad x_3=\bar\vartheta h.\]

\begin{lemma}\label{lem:casimir-h-nilpotence}
Let $n$ be as above, and let $M$ be a finite $D_{n,\zeta}$-module.
If its continuous dual $M'$ has an infinitesimal character, then
\[
 h^{2p^m}\mathsf S(M)=0
\]
for every choice of good filtration and integral graded lattice.
\end{lemma}
\begin{proof}
We may assume $M\ne0$. By direct computation, we get the Lie bracket relations
\[
 [X_1,X_2]=2X_3,\qquad
 [X_1,X_3]=2pX_2,\qquad
 [X_2,X_3]=2paX_1.
\]
Hence the Casimir element
\[
 \Omega=X_3^2+paX_1^2-pX_2^2
\]
is central in the enveloping algebra. The anti-involution defining
the dual action sends $X_i$ to $-X_i$ and fixes $\Omega$.
Thus the infinitesimal character on $M'$ gives
$(\Omega-\lambda)M=0$ for some $\lambda\in E$.

Since $ne=p^m$, the logarithm computation in Lemma~\ref{lem:cofinal-logarithm-symbol} gives us
\[
 \deg_n(X_3^2)=\deg_n(paX_1^2)=\deg_n(pX_2^2)
 =2p^m-2ne m.
\]
Multiplying the initial form of $\Omega$ by $\sigma_n(p)^{2m}$ gives
\[
 \sigma_n(p)^{2m}\sigma_n(\Omega)
 =x_3^{2p^m}+\sigma_n(p)
       \bigl(\bar a\,x_1^{2p^m}-x_2^{2p^m}\bigr).
\]
Modulo $\epsilon$, this is $x_3^{2p^m}$, since
$\sigma_n(p)\in k^\times\epsilon^e$.

We now consider $\sigma_n(\Omega-\lambda)$, $\lambda\in E$.
We have $\deg_n(\lambda)=nev_p(\lambda)$.

If $\deg_n(\lambda)<\deg_n(\Omega)$,
then $\sigma_n(\Omega-\lambda)=-\sigma_n(\lambda)\in k^\times\cdot\epsilon^{\mathbb{Z}}$ annihilates $\mathsf{S}(M)$.
Hence $\mathsf{S}(M)=0$ and therefore $M=0$, a contradiction.
Thus $\sigma_n(p^{2m}(\Omega-\lambda))\ \mathrm{mod}\ \epsilon=\sigma_n(p^{2m}\Omega)=x_3^{2p^m}$ annihilates $\mathsf{S}(M)$.
Finally, the relation $x_3\in k^\times h$ proves the assertion.
\end{proof}

We can now use the characteristic cycle over
$S/(h)=k[y,z]$ by applying Subsection~\ref{sec:uniform} with the
commutator ideal $(h)$. For a module $M$ as in the preceding
lemma, put $N=\mathsf S(M)$. If $M$ has canonical dimension one,
then
\[
 \dim\Supp_{k[y,z]}(\gr_{(h)}N)
 =3-j_S(N)=3-j_{D_{n,\zeta}}(M)=1.
\]
Here we use Lemma~\ref{lem:residual-grade} with $I=(h)$ and
Corollary~\ref{cor:quaternionic-grade}.

Write the one-dimensional cycle from Theorem~\ref{thm:additivity},
with the ideal $(h)$, as
\[
 \CC_1(M)=\sum_f m_f[V(f)],
\]
where the $f$ are irreducible homogeneous polynomials in $k[y,z]$.
We define the \emph{geometric multiplicity} of the cycle $\mathrm{CC}_1(M)$ by
\[
 \mu(M)\coloneqq\deg\CC_1(M)=\sum_f m_f\deg f,
\]
where we use the standard polynomial degree $\deg y=\deg z=1$.
This is the Hilbert multiplicity of
$\gr_{(h)}N$ over $k[y,z]$, and Theorem~\ref{thm:main} shows that it is independent
of the good filtration and integral graded lattice at this
fixed radius. We use the same convention for a finite graded
$S$-module annihilated by a power of $h$.

We now use equivariance to determine one-dimensional cycles of small degree.

\begin{lemma}\label{lem:geometric-two-J-support}
Let $N$ be a finite graded $S$-module annihilated by a power of
$h$, with a compatible $\widetilde H$-action. Suppose
$\gr_{(h)}N$ has support dimension 1 and geometric multiplicity $\mu(N)<\frac{p+1}{2}$.
Then $N$ is annihilated by a power of $J$, and
\[
 \CC_J(N)=\frac{\mu(N)}{2}([V(y)]+[V(z)]),\qquad \mu(N)<\frac{p+1}{2}.
\]
\end{lemma}
\begin{proof}
Let $q=\prod_f f^{m_f}$ be the homogeneous polynomial defining
the one-dimensional cycle of $\gr_{(h)}N$. The polynomial degree of $q$ is $\mu(N)$.
The cycle is $\widetilde H$-stable, so the line $kq$ is stable
under $H$. The monomials $y^{\mu(N)}, y^{\mu(N)-1}z,\dots,(yz)^{\mu(N)/2},\dots,z^{\mu(N)}$ have respective characters
\[
 \alpha^{\mu(N)}, \alpha^{\mu(N)-2}, \dots, \alpha^2,1,\alpha^{-2}, \dots, \alpha^{-\mu(N)}.
\]
These are distinct, since $\alpha$ has order $p+1$ and $p\geq5$.
Thus $q$ is one of these monomials. Conjugation by
$\varpi_D$ interchanges $y,z$ and preserves the cycle, so
$q=(yz)^{\mu(N)/2}$.

Any zero-dimensional homogeneous support is the origin. Hence
\[
 \Supp_{k[y,z]}\gr_{(h)}N\subseteq V(yz).
\]
Choose integers $a,b\geq1$ such that
$h^aN=0$ and $(yz)^b\gr_{(h)}N=0$. Since $J=(h,yz)$,
\[
 J^b h^iN\subseteq h^{i+1}N,
 \qquad J^{ab}N=0.
\]
We refine the $h$-filtration by powers of $yz$. The resulting
factors are killed by $J$, so we apply Lemma~\ref{lem:local-length}
to compare their generic lengths with those from the $J$-filtration.

The coefficients of the two axes are therefore the same as in the
original cycle, so both equal $\mu(N)/2$. This proves the asserted cycle and
multiplicity, with the convention of
\cite[\S3.2]{HuWangPadic}.
\end{proof}

We combine these results to get the following criterion.

\begin{proposition}\label{prop:infinitesimal-geometric-two}
Let $n=p^m/e$ with $m$ sufficiently large as above.
Let $M$ be a finite $\mathscr D_{n,\zeta}$-module 
with canonical dimension 1 and geometric
multiplicity 2. Suppose that
$M'$ has an infinitesimal character. Then every residual symbol $\mathsf S(M)$ is
annihilated by a power of $J$. Moreover,
\[
 \CC_J(M)=[V(y)]+[V(z)],
 \qquad\mu(M)=2.
\]
\end{proposition}
\begin{proof}
By Lemma~\ref{lem:casimir-h-nilpotence}, the residual symbol is
annihilated by a power of $h$. Choose a $\widetilde H$-stable
good filtration and integral graded lattice, as in
Lemma~\ref{lem:equivariant-module-filtration}. The $h$-adic filtration of the resulting residual symbol
has one-dimensional support and geometric multiplicity two.

We apply Lemma~\ref{lem:geometric-two-J-support}. Finally,
by Corollary~\ref{thm:quaternionic-CC}, the conclusion is independent of the
choice of filtration and lattice.
\end{proof}

Together with Theorem~\ref{thm:no-two}, this excludes geometric
multiplicity two along the cofinal sequence.

\begin{corollary}\label{cor:no-geometric-two}
There is no admissible locally analytic representation $V$ of
$D^\times$ with central character and infinitesimal character
such that, for all sufficiently large $n=p^m/e$,
the stage $(V'_b)_n$ has canonical dimension 1 and geometric
multiplicity 2.
\end{corollary}
\begin{proof}
Choose $m$ sufficiently large that
Proposition~\ref{prop:infinitesimal-geometric-two} applies and
the bound in Theorem~\ref{thm:no-two} holds. The proposition gives us
\[
 J^a\mathsf S((V'_b)_n)=0\quad\text{for some }a,
 \qquad \mu(\mathsf S((V'_b)_n))=2.
\]
This contradicts Theorem~\ref{thm:no-two}.
\end{proof}

\section{Application to Banach representations and the Jacquet-Langlands functor}
\label{sec:application}
In this section, we first give an irreducibility criterion for locally analytic
vectors in Banach representations, and then apply it to
$\Tau(\rho)$ using the comparison between $\JL(\rho)$ and
$\JL(\bar\rho)$.

\subsection{An irreducibility criterion for analytic vectors in Banach representations}
\label{sec:banach-criterion}
We first recall the required consequence of Sorgdrager's result.

\begin{proposition}\label{prop:sorgdrager-banach-support}
Let $\Theta\subset\Pi$ be a $D^\times$-stable open bounded lattice
in an admissible Banach representation with central character.
Suppose that $\Pi^{\mathrm{an}}$ has an infinitesimal character.
For any good $\omega$-filtration $F$, put
\[
 N=\gr_F((\Theta/\varpi_E\Theta)^\vee).
\]
Then
\[
 J^aN=0\text{ for some }a,\qquad
 \dim\Supp_{S/J}(\gr_JN)\leq1,
\]
with equality if $\Pi$ is infinite-dimensional.
\end{proposition}
\begin{proof}
We refer to the computation in \cite[\S6.6, Lemma~6.3; \S6.7, Theorem~6.2, equation (22)]{Sorgdrager2026} in the quaternion case over $\mathbb{Q}_p$. 
Using the notations in the previous section~\ref{sec:infinitesimal-geometric-two} before lemma~\ref{lem:casimir-h-nilpotence}, we know that for sufficiently large $n=p^m/e$ and $a=\vartheta^2\in\mathbb{Z}_p^\times$, 
\[N\text{ is annihilated by } (x_3^{2p^m}, \bar{a}x_1^{2p^m}-x_2^{2p^m}).\]
Sorgdrager uses the symbols $\mathbf w_1$ and $\mathbf{w}_2$ in this case, since $e_K=f_K=1$ for $K=\mathbb{Q}_p$; see \cite{Sorgdrager2026}.
By a direct computation in the coordinates of $S$, we get that $N$ is annihilated by $(h^{2p^m}, y^{p^m}z^{p^m})$, as $\bar{a}x_1^{2p^m}-x_2^{2p^m}=4\bar{a}y^{p^m}z^{p^m}$.
After a finite $(h)$ grading, we then deduce that $N$ is annihilated by a power of $J$. 

The dimension assertion follows from the fact that $S/J=k[y,z]/(yz)$ is of dimension 1 and the corollary~\ref{cor:D-zero}.
\end{proof}

We next choose a lattice with a pure associated graded module.

\begin{lemma}\label{lem:iwasawa-pure-lattice}
Let $V$ be a finite $\Lambda_{D,\zeta}$-module, pure of grade $j$, 
with a finite $\Lambda_{D, \zeta}^\circ$-lattice.
Then there is a commensurable $\Lambda_{D, \zeta}^\circ$-lattice $V^\circ$ such that
$V^\circ/\varpi_EV^\circ$ is pure of grade $j$ over $\overline\Lambda_{D, \bar\zeta}$
and admits a good $\omega$-filtration whose
associated graded module is pure of grade $j$ over $S\rtimes_{\bar{\zeta}}\Delta$.
\end{lemma}
\begin{proof}
We apply Theorem~\ref{thm:gabber-filtration} to the $\varpi_E$-adic filtration
of the given $\Lambda_{D, \zeta}^\circ$-lattice.
The resulting Gabber filtration $G$ is a good $\Lambda_{D, \zeta}^\circ$-filtration and has a pure associated
graded module. By definition it satisfies $G^iV=\varpi_E^iG^0V$. We set $V^\circ=G^0V$. Taking graded pieces along $G$, we have
\[
 \gr_GV=(V^\circ/\varpi_EV^\circ)[\epsilon^{\pm1}].
\]
For every finite $\overline\Lambda_{\bar\zeta}$-module $U$,
by the flat base change along $\overline\Lambda_{\bar\zeta}\to\overline\Lambda_{\bar\zeta}[\epsilon^{\pm1}]=\gr_{\varpi_E}\Lambda_{\zeta}$, we have
\[
 \operatorname{Ext}_{\overline\Lambda_{\bar\zeta}[\epsilon^{\pm1}]}^i
 \bigl(U[\epsilon^{\pm1}],
       \overline\Lambda_{\bar\zeta}[\epsilon^{\pm1}]\bigr)
 \cong\operatorname{Ext}_{\overline\Lambda_{\bar\zeta}}^i
 \bigl(U,\overline\Lambda_{\bar\zeta}\bigr)[\epsilon^{\pm1}].
\]
Since $\gr_GV$ is pure of grade $j$, we have
\[
 0\ne U\subseteq V^\circ/\varpi_EV^\circ
 \quad\Rightarrow\quad
 j_{\overline\Lambda_{\bar\zeta}}(U)
 =j_{\overline\Lambda_{\bar\zeta}[\epsilon^{\pm1}]}
      \bigl(U[\epsilon^{\pm1}]\bigr)=j.
\]
Hence $V^\circ/\varpi_EV^\circ$ is pure of grade $j$.
We then apply Theorem~\ref{thm:gabber-filtration} 
to the complete $\omega$-filtration on $\overline\Lambda_{D, \bar\zeta}$. 
Then we get the Gabber filtration over $\overline{\Lambda}_{D, \bar\zeta}$,
which is stable under $\widetilde{H}$ and has a pure associated graded module. 
\end{proof}

For $\Pi$ and $\Theta$ as in
Proposition~\ref{prop:sorgdrager-banach-support}, we put
$\overline M=(\Theta/\varpi_E\Theta)^\vee$ and define the Hilbert
multiplicity by
\[
 \dim_k(\overline M/\mD^r\overline M)
 =\mu(\Pi)r+O(1)\qquad(r\longrightarrow\infty).
\]
Thus $\mu(\Pi)=0$ when $\Pi$ is finite-dimensional.

The coefficient $\mu(\Pi)$ is additive and independent of the choice
of good filtration and Banach lattice, by the same lattice
comparison as in Proposition~\ref{prop:lattice-invariance}.
By Proposition~\ref{prop:sorgdrager-banach-support},
a finite $J$-filtration then gives us
\[
 \mu(\Pi)=\deg\CC_J(\gr_F\overline M)=\mu(\gr_F\overline M),
\]
in agreement with \cite[\S3.2 and Lemma~3.11]{HuWangPadic}. Here the degree
of $J$ is taken in the $\mathfrak{m}$-adic sense, i.e. $\deg y=\deg z=1$.

We can now state the criterion using the conventions above.

\begin{corollary}\label{cor:banach-four-analytic-irreducible}
Let $\Pi$ be an admissible unitary Banach representation of
$D^\times$ with central character $\zeta^{-1}$.
Suppose that $\Pi'$ is pure over $\Lambda_\zeta$ and
$\Pi^{\mathrm{an}}$ has an infinitesimal character. If
\[
 \Pi^{1\text{-}\mathrm{alg}}=0,\qquad \mu(\Pi)=4,
\]
then $\Pi^{\mathrm{an}}$ is topologically irreducible.
\end{corollary}
\begin{proof}
By Lemma~\ref{lem:iwasawa-pure-lattice} and Schikhof duality,
we choose a $D^\times$-stable compact dual lattice $M^\circ\subset\Pi'$ and a good
$\omega$-filtration $F$ with a pure associated graded module.

Proposition~\ref{prop:sorgdrager-banach-support} tells us that $\gr_F(M^\circ/\varpi_EM^\circ)$ is annihilated by a power of $J$
and has canonical dimension 1 with $\mu(M)=4$. We lift the homogeneous generators of $\gr_F(M^\circ/\varpi_EM^\circ)$, with their fixed degrees, to filtered generators in $M^\circ$, say $m_i$.

We put $M=(\Pi^{\mathrm{an}})'_b$. By \cite[Theorem~7.1]{ST03}, we have
$M_n=D_{n,\zeta}\otimes_{\Lambda_\zeta^\circ}M^\circ$.
Then Lemma~\ref{lem:large-radius-recovery} gives us with the good filtration on $M_n$ induced by the generators $1\otimes m_i$, 
for all large $n$,
\[
 \mathsf S(M_n)\cong\gr_F(M^\circ/\varpi_EM^\circ),\qquad
 \mu(M_n)=4.
\]
Thus these symbols are pure of grade 2 and killed by a power of $J$.

Suppose that $0\ne L\subsetneq M$ is a closed $D^\times$-stable
submodule, and let $L_n,Q_n$ be the stages of $L,M/L$.
For all large $n$, both are nonzero. By Lemma~\ref{lem:two-stage-strictness},
the induced subspace and quotient filtrations give the exact sequence
\[
 0\longrightarrow\mathsf S(L_n)
 \longrightarrow\mathsf S(M_n)
 \longrightarrow\mathsf S(Q_n)\longrightarrow0.
\]
Since $\mathsf{S}(M_n)$ is pure, Corollary~\ref{cor:D-zero} gives us
$\mu(L_n)>0$.
If $\mu(Q_n)=0$, the same corollary shows that $Q_n$
is finite-dimensional. The dense map $M\to Q_n$ is then surjective,
so
\[
 0\ne Q_n'\hookrightarrow\Pi^{\mathrm{an}}\hookrightarrow\Pi.
\]
This contradicts $\Pi^{1\text{-}\mathrm{alg}}=0$, since finite-dimensional
$D^\times$-representations are locally algebraic on $D^{\times,1}$
\cite[Lemma~4.22]{LiuWang2025}.
By additivity and Corollary~\ref{cor:D-symmetry}, we have
\[
 \mu(L_n)=\mu(Q_n)=2.
\]
The symbols $\mathsf S(Q_n)$ are killed by a power of $J$ with $\mu(Q_n)=2$ 
for all sufficiently large $n$, contradicting Corollary~\ref{cor:no-eventual-two}.
The assertion follows by locally analytic duality \cite[Theorem~6.3]{ST03}.
\end{proof}

\begin{remark}
The condition $\Pi^{1\text{-}\mathrm{alg}}=0$ excludes
finite-dimensional subrepresentations even for a nonalgebraic
central character; see the correction preceding
\cite[Theorem~5.3]{LiuWang2025}.
\end{remark}

\subsection{The Jacquet-Langlands functor}
\label{sec:HW-tori}
We first compare the rational $D^\times$-representation
$\JL(\rho)$ with the residual representation $\JL(\bar\rho)$.

\paragraph{Assumptions on \texorpdfstring{$\rho$}{rho}}
We first state the assumptions on $\rho$ as follows:
\begin{hypothesis}\label{hypothesis}
Let $\rho\colon\mathrm{Gal}_{\mathbb{Q}_p}\to\mathrm{GL}_2(E)$ be a two dimensional absolutely irreducible continuous representation. Suppose that 
there exist a $\mathrm{Gal}_{\mathbb{Q}_p}$-stable lattice $T_\rho\subset\rho$ with $\bar{\rho}=T_\rho/\varpi_E$ satisfies the one of the following genericity conditions:
\begin{itemize}
      \item[Type~(I)] 
             $\mathscr B(\bar\rho)=\{\pi(\bar{\rho})\}$ is of block type~(I), i.e. the supersingular block, which corresponds to $\bar{\rho}$ being irreducible;  $\bar{\rho}$ is generic in the sense $\bar{\rho}|_{I_{\mathbb{Q}_p}}\sim\begin{psmallmatrix}
                  \omega_2^{r+1} & 0 \\ 0 & \omega_2^{p(r+1)}
            \end{psmallmatrix}\otimes\chi$ for $2\le r\le p-3$, where $\omega_2$ is Serre's fundamental character of niveau $2$.
            
            \item[Type~(II)] $\mathscr B(\bar\rho)=\{\mathrm{Ind}_B^{\mathrm{GL_2}}(\chi_1\otimes\chi_2\bar{\chi}_{\mathrm{cycl}}^{-1}), \mathrm{Ind}_B^{\mathrm{GL}_2}(\chi_2\otimes\chi_1\bar{\chi}_{\mathrm{cycl}}^{-1})\}$ 
                  with $\chi_1\chi_2^{-1}\neq1,\bar{\chi}_{\mathrm{cycl}}^{\pm1}$ is of block type~(II). And $\bar{\rho}$ is generic in the sense $\bar{\rho}|_{I_{\mathbb{Q}_p}}\sim\begin{psmallmatrix}
                  \bar{\chi}_{\mathrm{cycl}}^{r+1} & * \\ 0 & 1
            \end{psmallmatrix}\otimes\chi$ for $1\le r\le p-3$ and $\bar{\rho}^{\mathrm{ss}}\neq(\bar{\chi}_{\mathrm{cycl}}\oplus1)\otimes\chi$.
\end{itemize}

\end{hypothesis}

Let $\mathcal S^i$ denote Scholze's functors \cite{Scholze2018},
\[
 \mathcal S^i(\pi)
 :=\mathrm H^i_{\mathrm{\acute et}}
   (\mathbb P^1_{\mathbb C_p},\mathcal F_\pi).
\]
We use Hu--Wang's residual representation $\JL(\bar\rho)$, with the
split convention of \cite[Remark~8.36]{HuWangModP}. More explicitly,
\[
 \JL(\bar\rho):=
 \begin{cases}
 \Hom_{\mathrm{Gal}_{\mathbb Q_p}}
   (\bar\rho(-1),\mathcal S^1(\pi(\bar\rho))),
   &\operatorname{End}_{\mathrm{Gal}_{\mathbb Q_p}}(\bar\rho)=k,\\[2pt]
 \Hom_{\mathrm{Gal}_{\mathbb Q_p}}
   (\chi_1(-1),\mathcal S^1(\pi(\bar\rho))),
   &\bar\rho=\chi_1\oplus\chi_2.
 \end{cases}
\]
By \cite[Theorem~4.20]{HuWangPadic},
\[
 \mathcal S^1(\pi(\bar\rho))
 \cong\bar\rho(-1)\otimes_k\JL(\bar\rho).
\]
In the split case, either character gives an isomorphic
$D^\times$-representation, whereas
\[
 \Hom_{\mathrm{Gal}_{\mathbb Q_p}}
 (\bar\rho(-1),\mathcal S^1(\pi(\bar\rho)))
 \cong\operatorname{End}_{\mathrm{Gal}_{\mathbb Q_p}}(\bar\rho)
       \otimes_k\JL(\bar\rho)
 \cong\JL(\bar\rho)^{\oplus2}.
\]

In type~\textup{(II)}, write
$\bar\rho^{\mathrm{ss}}=\chi_1\oplus\chi_2$ and put
\[
 \pi_1:=\operatorname{Ind}_B^{\GL_2(\mathbb Q_p)}
    (\chi_2\otimes\chi_1\bar\chi_{\mathrm{cycl}}^{-1}),\qquad
 \pi_2:=\operatorname{Ind}_B^{\GL_2(\mathbb Q_p)}
    (\chi_1\otimes\chi_2\bar\chi_{\mathrm{cycl}}^{-1}).
\]
The proof of \cite[Theorem~8.12]{HuWangModP}, together with
\cite[Theorem~4.20]{HuWangPadic}, gives
\[
 \mathcal S^1(\pi_i)\cong\chi_i(-1)\otimes_k\JL(\bar\rho)
 \quad(i=1,2).
\]
These are the residual compatibilities used below.

Choose any $\GL_2(\mathbb Q_p)$-stable open bounded lattice
$\Theta_\rho\subset\Pi(\rho)$ and put
$\bar\pi_\rho:=\Theta_\rho/\varpi_E\Theta_\rho$.
Compatibility of the $p$-adic and mod-$p$ local Langlands
correspondences \cite[Proposition~III.55 and Remark~III.56]{ColmezDospinescu2014} gives
\[
 \bar\pi_\rho^{\mathrm{ss}}\cong
 \begin{cases}
  \pi(\bar\rho),&\text{in type~(I)},\\
  \pi_1\oplus\pi_2,&\text{in type~(II)}.
 \end{cases}
\]
By \cite[Theorem~4.11]{HuWangPadic} and the long exact sequence,
\[
 \mathcal S^0(\bar\pi_\rho)=\mathcal S^2(\bar\pi_\rho)=0.
\]
The integral and rational Scholze functors are therefore given by
\[
 \check{\mathcal S}^1(\Theta_\rho)
 :=\varprojlim_m\mathcal S^1(\Theta_\rho/\varpi_E^m\Theta_\rho),
\qquad
 \check{\mathcal S}^1(\Pi(\rho))
 =\check{\mathcal S}^1(\Theta_\rho)[1/\varpi_E].
\]
Then \cite[Proposition~4.14]{HuWangPadic} gives the exact sequence
\[
 0\longrightarrow\check{\mathcal S}^1(\Theta_\rho)
 \xrightarrow{\ \varpi_E\ }\check{\mathcal S}^1(\Theta_\rho)
 \longrightarrow\mathcal S^1(\bar\pi_\rho)\longrightarrow0.
\]
In our normalization, rational typicity \cite[Corollary~4.20]{LiuWang2025} gives us
\begin{equation}\label{eq:rational-S-typicity}
 \check{\mathcal S}^1(\Pi(\rho))\cong\rho(-1)\widehat\otimes_E\JL(\rho)\quad\Rightarrow\quad\JL(\rho)=\Hom_{\mathrm{Gal}_{\mathbb{Q}_p}}(\rho(-1), \check{\mathcal S}^1(\Pi(\rho))).
\end{equation}

\begin{proposition}\label{prop:integral-jl-lattice}
Assume that $\bar\rho$ is generic of type~\textup{(I)} or type~\textup{(II)}.  There exists an
open bounded $D^\times$-stable lattice
\[
 \JL(\rho)^\circ\subset\JL(\rho)
\]
such that
\[
 \JL(\rho)^\circ[1/\varpi_E]=\JL(\rho),
 \qquad
 \JL(\rho)^\circ/\varpi_E\JL(\rho)^\circ
 \cong\JL(\bar\rho).
\]
\end{proposition}

\begin{proof}
We first choose an element $g_0\in\mathrm{Gal}_{\mathbb Q_p}$ for
which $\bar\rho(g_0)$ has two distinct eigenvalues.  Such an element exists since $\bar{\rho}$ is either irreducible or an extension of two different characters. 
Twisting by $(-1)$ does not change this property.  Denote
the two eigenvalues of $\bar\rho(-1)(g_0)$ by
$\bar\lambda_1$ and $\bar\lambda_2$.

The characteristic polynomial of $\rho(-1)(g_0)$ therefore has two
distinct roots modulo $\varpi_E$.  After enlarging $E$ if necessary,
Hensel's lemma gives roots $\lambda_1,\lambda_2\in\cO_E$ satisfying
\[
 \lambda_i\bmod\varpi_E=\bar\lambda_i,
 \qquad
 \lambda_1-\lambda_2\in\cO_E^\times.
\]
Let $g_0$ also denote its action on $\check{\mathcal S}^1(\Theta_\rho)$, and set
\[
 \mathrm{pr}_1:=\frac{g_0-\lambda_2}{\lambda_1-\lambda_2}.
\]
The denominator is a unit, so $\mathrm{pr}_1$ preserves
$\check{\mathcal S}^1(\Theta_\rho)$.  By \eqref{eq:rational-S-typicity}, the
operator $g_0$ acts on $\check{\mathcal S}^1(\Pi(\rho))$ through
$\rho(-1)(g_0)\otimes 1$.  The Cayley--Hamilton identity thus gives
\[
 (g_0-\lambda_1)(g_0-\lambda_2)=0
 \quad\text{on }\check{\mathcal S}^1(\Pi(\rho)).
\]
Since $\check{\mathcal S}^1(\Theta_\rho)$ is $\varpi_E$-torsion-free, the same
identity holds on $\check{\mathcal S}^1(\Theta_\rho)$.  Consequently,
\[
 \mathrm{pr}_1^2=\mathrm{pr}_1.
\]
The actions of $\mathrm{Gal}_{\mathbb Q_p}$ and $D^\times$ commute,
so $\mathrm{pr}_1\check{\mathcal S}^1(\Theta_\rho)$ is $D^\times$-stable.

From \eqref{eq:rational-S-typicity} we obtain
\[
 \mathrm{pr}_1\check{\mathcal S}^1(\Pi(\rho))
 \cong
 \bigl(\mathrm{pr}_1\rho(-1)\bigr)\widehat\otimes_E\JL(\rho).
\]
The eigenspace $\mathrm{pr}_1\rho(-1)$ is one-dimensional.  After choosing an
identification $\mathrm{pr}_1\rho(-1)\cong E$, we obtain a $D^\times$-equivariant
isomorphism
\[
 \mathrm{pr}_1\check{\mathcal S}^1(\Pi(\rho))\cong\JL(\rho).
\]
We define \[\JL(\rho)^\circ\coloneqq\mathrm{Im}(\mathrm{pr}_1\check{\mathcal S}^1(\Theta_\rho)\to \mathrm{pr}_1\check{\mathcal S}^1(\Pi(\rho))\cong\JL(\rho)).\]
It is an open
bounded $D^\times$-stable lattice in $\JL(\rho)$, i.e. $\JL(\rho)^\circ[1/\varpi_E]=\JL(\rho)$.

Since $\mathrm{pr}_1$ is an integral idempotent, reduction commutes
with its image:
\[
 \JL(\rho)^\circ/\varpi_E\JL(\rho)^\circ
 \cong\bar{\mathrm{pr}}_1
       (\check{\mathcal S}^1(\Theta_\rho)/\varpi_E\check{\mathcal S}^1(\Theta_\rho))
 \cong\bar{\mathrm{pr}}_1\mathcal S^1(\bar\pi_\rho).
\]
In type~\textup{(I)}, $\bar\pi_\rho\cong\pi(\bar\rho)$, and
\cite[Theorem~4.20]{HuWangPadic} gives
\[
 \bar{\mathrm{pr}}_1\mathcal S^1(\bar\pi_\rho)
 \cong\bigl(\bar{\mathrm{pr}}_1\bar\rho(-1)\bigr)
       \otimes_k\JL(\bar\rho)
 \cong\JL(\bar\rho).
\]
In type~\textup{(II)}, there is an exact sequence
\[
 0\longrightarrow\pi_i\longrightarrow\bar\pi_\rho
 \longrightarrow\pi_j\longrightarrow0,
 \qquad\{i,j\}=\{1,2\}.
\]
The vanishing of $\mathcal S^0(\pi_j)$ and $\mathcal S^2(\pi_i)$,
together with the principal-series identities above, gives
\[
 0\longrightarrow\chi_i(-1)\otimes_k\JL(\bar\rho)
 \longrightarrow\mathcal S^1(\bar\pi_\rho)
 \longrightarrow\chi_j(-1)\otimes_k\JL(\bar\rho)
 \longrightarrow0.
\]
Label the eigenvalues so that
$\bar\lambda_1=\chi_1(-1)(g_0)$. Applying
$\bar{\mathrm{pr}}_1$ is exact; it acts as the identity on the
$\chi_1(-1)$-isotypic term and as zero on the $\chi_2(-1)$-isotypic term. Hence
\[
 \bar{\mathrm{pr}}_1\mathcal S^1(\bar\pi_\rho)
 \cong\JL(\bar\rho)
\]
for either order of the two factors, including when $\bar\rho$ is split.
All these isomorphisms are $D^\times$-equivariant, so their duals also
identify the $\mD$-adic filtrations.
\end{proof}

\paragraph{Residual symbols and multiplicity}

Let $\JL(\rho)^\circ$ be the lattice of
Proposition~\ref{prop:integral-jl-lattice}, and set
\[
 \JL(\rho)^{\circ,\mathbf d}
 :=\Hom_{\cO_E}^{\mathrm{cont}}(\JL(\rho)^\circ,\cO_E).
\]
This is a finite $\Lambda_{D,\zeta}^\circ$-module and spans the
continuous dual $\JL(\rho)'$ and it satisfies
\[
 \JL(\rho)^{\circ,\mathbf d}/\varpi_E
 \cong\JL(\bar\rho)^\vee.
\]
By Lemma~\ref{omega=m-adic}, the $\omega$-grading coincides with the
$\mD$-adic grading after rescaling. Hence
\cite[Theorem~4.19]{HuWangPadic} gives
\[
 \gr_\omega^*\bigl(\JL(\bar\rho)^\vee\bigr)
 \cong\bigoplus_{\chi\in W_D(\bar\rho)}
       \chi^\vee\otimes_k A_0/\mathfrak a(\chi),
 \qquad \deg_\omega y=\deg_\omega z=\tfrac12.
\]
This graded $S$-module is Cohen--Macaulay of grade two, hence
pure, by \cite[Theorem~4.19 and Lemma~6.2]{HuWangPadic}.
It is annihilated by $J$ and has multiplicity four.

The integral comparison and the definition of Banach multiplicity give us
\[
 \mu(\JL(\rho))=\mu(\JL(\bar\rho))=4.
\]

We also need purity over the rational Iwasawa algebra.

\begin{lemma}\label{lem:jl-iwasawa-pure}
Under Hypothesis~\ref{hypothesis}, the continuous dual
$\JL(\rho)'$ is pure of grade two over $\Lambda_\zeta$.
\end{lemma}
\begin{proof}
Let $W\subseteq\JL(\rho)'$ be a nonzero
$\Lambda_\zeta$-submodule, and put
\[
 W^\circ=W\cap\JL(\rho)^{\circ,\mathbf d}.
\]
It is finite over $\Lambda_\zeta^\circ$, spans $W$, and is
$\varpi_E$-saturated. Thus
\[
 0\ne W^\circ/\varpi_EW^\circ
 \hookrightarrow\JL(\bar\rho)^\vee.
\]
Give the left-hand module the induced good $\omega$-filtration.
Its associated graded embeds in the pure grade-two module
$\gr_\omega^*(\JL(\bar\rho)^\vee)$. The filtered grade comparison
\cite[Chapter~III, Theorem~2.5.2]{LvO96} gives
\[
 j_{\overline\Lambda_{\bar\zeta}}
       (W^\circ/\varpi_EW^\circ)=2.
\]

Now filter $W$ by $\varpi_E^iW^\circ$. As in the proof of
Lemma~\ref{lem:iwasawa-pure-lattice}, its associated graded is
$(W^\circ/\varpi_EW^\circ)[\epsilon^{\pm1}]$.
Laurent extension and the same filtered grade comparison give
\[
 j_{\Lambda_\zeta}(W)
 =j_{\overline\Lambda_{\bar\zeta}[\epsilon^{\pm1}]}
      \bigl((W^\circ/\varpi_EW^\circ)[\epsilon^{\pm1}]\bigr)
 =j_{\overline\Lambda_{\bar\zeta}}
      (W^\circ/\varpi_EW^\circ)=2.
\]
This proves purity.
\end{proof}

\paragraph{Irreducibility of locally analytic Jacquet-Langlands}
\label{sec:analytic-irreducibility}

We now apply the criterion in Corollary~\ref{cor:banach-four-analytic-irreducible} to $\Tau(\rho)$ and $\Tau(\rho)^{\mathrm{an}}$.

\begin{corollary}\label{thm:tau-analytic-irreducible}
Under Hypothesis~\ref{hypothesis}, the admissible
locally analytic representation $\Tau(\rho)^{\mathrm{an}}$ is
topologically irreducible.
\end{corollary}
\begin{proof}
Recall that
\[
 \Tau=\JL(\rho)/\JL(\rho)^{1\text{-}\mathrm{alg}},
 \qquad \dim_E\JL(\rho)^{1\text{-}\mathrm{alg}}<\infty
\]
by \cite[Theorem~4.23]{LiuWang2025}. The quotient is admissible
and unitary, with the same central character. Its continuous
dual is a submodule of $\JL(\rho)'$, so it is pure of grade two
by Lemma~\ref{lem:jl-iwasawa-pure}. The finite-dimensional term
has zero linear Hilbert coefficient. Additivity therefore gives
\[
 \mu(\Tau)=\mu(\JL(\rho))=4.
\]

We check that $\Tau^{1\text{-}\mathrm{alg}}=0$.
If $V\subseteq\Tau$ is a finite-dimensional
$D^\times$-subrepresentation, its inverse image $W$ fits into
\[
 0\longrightarrow\JL(\rho)^{1\text{-}\mathrm{alg}}
 \longrightarrow W\longrightarrow V\longrightarrow0.
\]
Then $W$ is finite-dimensional, hence locally algebraic for
$D^{\times,1}$ by \cite[Lemma~4.22]{LiuWang2025}.
Thus $W\subseteq\JL(\rho)^{1\text{-}\mathrm{alg}}$ and thus $V=0$.

By \cite[Corollary~8.16]{DPS23},
$\check{\mathcal S}^1(\Pi(\rho))^{\mathrm{an}}$ has an
infinitesimal character. The decomposition
\eqref{eq:rational-S-typicity} gives the same character on
$\JL(\rho)^{\mathrm{an}}$ and then passes it to $\Tau^{\mathrm{an}}$.
Now all the hypotheses of
Corollary~\ref{cor:banach-four-analytic-irreducible} are satisfied,
so $\Tau^{\mathrm{an}}$ is topologically irreducible.
\end{proof}

\printbibliography[heading=bibintoc]
\end{document}